\documentclass[90pt]{article}

\usepackage{verbatim}
\usepackage{amsmath}
\usepackage{bbm}
\usepackage{amsfonts}
\usepackage{amsthm}

\newtheorem{rem1}{Remark}[section]
\newtheorem{lem1}{Lemma}[section]
\newtheorem{cor1}{Corollary}[section]

\newtheorem{prop1}{Proposition}[section]
\newtheorem{thm1}{Theorem}[section]

\begin{document}
\title{Asymptotics of spectral quantities of Zakharov--Shabat operators}

\author{T. Kappeler\footnote{Supported in part by the Swiss National Science Foundation}, 
B. Schaad\footnote{Supported in part by the Swiss National Science Foundation}, 
P. Topalov\footnote{Supported in part by NSF DMS-0901443}}
\maketitle

\begin{abstract}
In this paper we provide new asymptotic estimates of various spectral quantities of Zakharov-Shabat operators 
on the circle. These estimates are uniform on bounded subsets of potentials in Sobolev spaces.
\end{abstract}

\section{Introduction}\label{Introduction}
In this paper we prove asymptotic estimates of various spectral quantities of Zakharov--Shabat (ZS) operators 
\[L(\varphi) = i \begin{pmatrix}
1&0\\0&-1
\end{pmatrix}\partial_x + \begin{pmatrix}
0& \varphi_1 \\
\varphi_2 &0
\end{pmatrix}\]
in one space dimension. These operators appear in the Lax pair formulation of the focusing and defocusing NLS equation and hence their 
spectral properties are relevant for the study of these equations. 
We assume that $\varphi= (\varphi_1,\varphi_2)$ is in $H^N_c= H^N\times H^N,\, N\in \mathbb{Z}_{\geq 0}$, 
where $H^N$ denotes the Sobolev space of 1-periodic complex-valued functions supplied with the standard Sobolev norm
$\|u\|_{H^N}:=\big(\sum_{j=0}^N\|\partial_x^j u\|_{L^2}^2\big)^{1/2}$, $\|u\|_{L^2}:=\int_0^1|u(x)|^2\,dx$. 
For a given potential $\varphi\in H^0_c\equiv L^2_c$, consider the operator $L(\varphi)$ with {\em periodic} boundary conditions on 
the interval $[0,2].$  Note that unless $\varphi_2=\overline{\varphi_1}$, $L(\varphi)$ is {\em not}  formally selfadjoint with respect to 
the $L^2$-inner product on $[0,2]$,
\[
\langle F, G\rangle = \frac{1}{2}\int_0^2 \big(F_1\overline{G_1}+ F_2\overline{G_2}\big)\,dx,
\]
where $F=(F_1, F_2)$ and $G=(G_1, G_2)$ are complex-valued $L^2$-functions on $[0,2]$.
In addition, we will also consider $L(\varphi)$ with {\em Dirichlet} boundary conditions on $[0,1]$ whose domain consists of 
all functions $F= (F_1, F_2)$ in $H^1([0,1],\mathbb{C})\times H^1([0,1],\mathbb{C})$ such that
\[ 
F_1(0)= F_2(0), \;F_1(1)= F_2(1)\,.
\]
The corresponding spectra, referred to as periodic, respectively Dirichlet spectrum of $L(\varphi),$
are discrete. The eigenvalues can be listed (with their algebraic multiplicities) as sequences of complex numbers
\[\cdots \preceq\lambda^-_{n}\preceq \lambda^+_n\preceq \lambda^-_{n+1} 
\preceq\lambda_{n+1}^+ \preceq \cdots \quad\text{and}\quad
\cdots \preceq \mu_n\preceq \mu_{n+1} \preceq \cdots\]
in lexicographic order $\preceq$ in such a way that  
\begin{align}\label{1.2}
\mu_n,\lambda_n^\pm= n\pi +\ell^2_n\quad \text{as}\;\;\; |n|\to \infty
\end{align} 
-- see e.g. \cite{GK}, Proposition 5.3 and Proposition 6.7.
Two complex numbers $a$ and  $b$ are {\em lexicographically} ordered $a\preceq b$,  
if $[\operatorname{Re}(a)< \operatorname{Re}(b)]$ or [$\operatorname{Re}(a)=\operatorname{Re}(b)$
and $[\operatorname{Im}(a)\leq\operatorname{Im}(b)$].  
The notation $\mu_n= n\pi +\ell^2_n$ means that $(\mu_n-n\pi)_{n\in \mathbb{Z}}$ is an $\ell^2$--sequence. 
Furthermore denote by $M(x,\lambda)\equiv M(x,\lambda,\varphi)$ the fundamental solution  
\[
M(x,\lambda)= \begin{pmatrix}
m_1(x,\lambda)& m_2(x,\lambda) \\
m_3(x,\lambda) & m_4(x,\lambda)
\end{pmatrix},\quad M(0,\lambda)= \begin{pmatrix}
1&0\\0&1
\end{pmatrix}, 
\]
of the linear system $L(\varphi) M=\lambda M$, $\lambda\in\mathbb{C}$.
For any $x\in \mathbb{R}$, $M(x,\lambda)$ is an entire function in $\lambda$. 
Let $\Delta(\lambda) \,[\delta(\lambda)]$ be the the trace [anti--trace] of $M(1,\lambda)$
\[
\Delta(\lambda):=m_1(1,\lambda)+m_4(1,\lambda),\quad \delta(\lambda):=m_2(1,\lambda)+m_3(1,\lambda)
\]
and set ${\dot\Delta}(\lambda):= \partial_\lambda \Delta(\lambda).$ 
The zeros of ${\dot\Delta}$ can be listed (with their multiplicities) as a sequence  
of complex numbers $\cdots \preceq \dot\lambda_n\preceq \dot\lambda_{n+1} \preceq \cdots$ in 
lexicographic order so that 
\begin{align}\label{1.4}
{\dot\lambda_n}= n\pi +\ell^2_n\quad \text{as}\; |n|\to \infty
\end{align} 
-- see e.g. \cite{GK}, Lemma 6.5.
Furthermore, let
$\tau_n:= (\lambda_n^++ \lambda_n^-) / 2$ and $\gamma_n:=\lambda_n^+ -\lambda_n^-.$
Note that by \eqref{1.2},
\begin{align}
\label{1.6} \tau_n= n\pi+ \ell^2_n \quad \text{and} \quad \gamma_n^2= \ell^1_n.
\end{align}

\noindent The  aim of this paper is to establish refined asymptotics of  $\mu_n$, $\lambda_n^{\pm}$, $\gamma_n^2$, $\tau_n$,  
and $\dot\lambda_n$ as $|n|\to \infty$ as well as asymptotics of other spectral quantities such as  $\Delta(\mu_n) $ and 
$\delta(\mu_n)$ for potentials in $H^N_c$ with  $N\in \mathbb{Z}_{\geq 1}.$ 
For any $s\geq 0$, consider the real subspace of $H^s_c$,
\[
H^s_r:=\big\{(u,\bar{u})\,|\,u\in H^s\big\}\,.
\]
For $\varphi\in H^0_r\equiv L^2_r$ the operator $L(\varphi)$ considered with periodic and Dirichlet boundary conditions as 
discussed above is selfadjoint.
In particular, all the quantities $\mu_n$, $\lambda_n^{\pm}$, $\tau_n$, and $\dot\lambda_n$ are real-valued. 
Denote by $\hat u(n)$, $n\in\mathbb{Z}$, the $n$-th Fourier coefficient
of a 1-periodic function $u\in H^0\equiv L^2$,  ${\hat u}_n:=\int_0^1 u(x)e^{-2\pi i n x}\,dx$.
\begin{thm1}\label{Theorem1.1} 
For $\varphi \in H^N_c$ with $N\geq 1,$  
\[
\mu_n= n\pi +\sum_{k=1}^{N+1}\frac{c_k}{n^k} +
\frac{1}{2}\big(\hat \varphi_1(-n) +\hat\varphi_2(n)\big) + \frac{\ell^2_n}{n^{N+1}} \quad \text{ as }\,\,
|n|\to \infty 
\]
uniformly on bounded sets of $H^N_c.$ The coefficients $c_k\equiv c_k(\varphi)$ are independent of the choice of $n$ and $N$
and can be represented as integrals of polynomials of $\varphi_1,\varphi_2$ and their derivatives up to order $k-1.$
\end{thm1}
\begin{rem1}
The coefficients $c_k$ can be computed inductively -- see Remark \ref{rem1}. 
One has $c_1 = \frac{1}{2\pi}\int_0^1 \varphi_1(t)\varphi_2(t) dt$ and 
$c_2 = \frac{i}{4\pi^2}\int_0^1 \varphi_1(t)\varphi_2'(t) dt.$
\end{rem1}
\begin{thm1}\label{Theorem1.2}
\begin{itemize}
\item[(i)] For $\varphi \in H^N_c$ with $N\in  \mathbb{Z}_{\geq 1}$,
\[
\{\lambda^+_n,\lambda^-_n\}= \Big\{ 
n\pi +\sum_{k=1}^{N+1}\frac{c_k}{n^k} \pm\sqrt{\hat \varphi_1(-n) \hat\varphi_2(n)} +\frac{\ell^2_n}{n^{N+\frac{1}{2}}}
\Big\}\quad \text{as}\,\,|n|\to \infty\,.
\]
uniformly on bounded sets of $H^N_c.$ 
\item[(ii)]
For $\varphi \in H^N_r$ with $N\in  \mathbb{Z}_{\geq 1},$ 
\[\lambda_n^\pm=n\pi +\sum_{k=1}^{N+1}\frac{c_k}{n^k} \pm\sqrt[+]{\hat \varphi_1(-n) \hat\varphi_2(n)} +
\frac{\ell^4_n}{n^{N+1}} \quad \text{as} \; |n|\to \infty\]
uniformly on bounded sets of $H^N_r.$
\end{itemize}
\noindent The coefficients $c_k$ are the same as in Theorem \ref{Theorem1.1}.

\end{thm1}
\begin{rem1}\label{Remark1}
Note that $\lambda^-_{n}\preceq \lambda^+_n$ whereas the two values of the square root $\sqrt{\hat \varphi_1(-n) \hat\varphi_2(n) }$ are 
not lexicographically ordered in a canonical way. For this reason, in item (i), the asymptotics of $\lambda_n^\pm$ are stated in terms of an 
equality of sets. In contrast, for $\varphi\in H^N_r,\,\hat \varphi_1(-n)=\overline{\hat \varphi}_2(n)$ and hence 
$\sqrt[+]{\hat \varphi_1(-n) \hat\varphi_2(n)}\geq 0$, allowing to specify the asymptotics as in (ii).
\end{rem1}
\noindent
As an immediate application of Theorem \ref{Theorem1.2} one gets the following
\begin{cor1}\label{Corollary1.3}
\begin{itemize}
\item[(i)] For $\varphi\in H^N_c$ with $N\in \mathbb{Z}_{\geq 1},$ 
\[\gamma_n= 2 \sqrt{\hat \varphi_1(-n) \hat\varphi_2(n)}+\frac{\ell^2_n}{n^{N+\frac{1}{2}}}\quad \text{as} \; |n|\to \infty \]
with the appropriate choice of the square root. The asymptotics hold
uniformly on bounded sets of $H^N_c.$ 
\item[(ii)] For $\varphi\in H^N_r$ with $N\in \mathbb{Z}_{\geq 1},$
\[
0\leq \gamma_n= 2|\varphi_1(-n)| +\frac{\ell^4_n}{n^{N+1}}\quad \text{as} \; |n|\to \infty 
\]
uniformly on bounded sets of $H^N_r.$
\end{itemize}
\end{cor1}

\medskip

\noindent In \cite{KST2} we need asymptotic estimates for $\tau_n=(\lambda_n^++\lambda_n^-) / 2.$ 
But the ones obtained from Theorem \ref{Theorem1.2} are not sufficient for our purposes.  
We derive the sharper estimates from asymptotic estimates of the zeros $(\dot\lambda_n)_{n\in\mathbb{Z}}$ of $\dot\Delta(\lambda).$
\begin{thm1}\label{Theorem1.4}
For $\varphi \in H^N_c$ with $N\geq 1,$ 
\begin{itemize}
\item[(i)] $ \dot\lambda_n=n\pi +\sum_{k=1}^{N+1}\frac{c_k}{n^k}+\frac{\ell^2_n}{n^{N+1}} \quad \text{as} \; |n|\to \infty$
\item[(ii)] $ \tau_n=n\pi +\sum_{k=1}^{N+1}\frac{c_k}{n^k}+\frac{\ell^2_n}{n^{N+1}} \quad \text{as} \; |n|\to \infty$
\end{itemize}
uniformly on bounded sets of $H^N_c.$
The coefficients $c_k$ are the same as in Theorem \ref{Theorem1.1}.
\end{thm1}
\noindent
Finally, in \cite{KST2} we also need asymptotic estimates for $ \Delta(\mu_n)$ and $ \delta(\mu_n)$.
\begin{thm1}\label{Theorem1.5}
For $\varphi \in H^N_c$ with $N\geq 1,$ 
\begin{itemize}
\item[(i)] $ \Delta(\mu_n)=(-1)^n2+ +\frac{\ell^2_n}{n^{N+1}} \quad \text{as} \; |n|\to \infty$
\item[(ii)] $ \delta(\mu_n)=(-1)^n i\,\big(\hat \varphi_1(-n)-\varphi_2(n) \big) +\frac{\ell^2_n}{n^{N+1}} \quad \text{as} \; |n|\to \infty$
\end{itemize}
uniformly on bounded sets of $H^N_c.$
\end{thm1}
\noindent
To prove the stated asymptotic estimates we need to define and study special solutions of 
$L(\varphi)F=\lambda F$ for $\lambda \in \mathbb{C}$ sufficiently large which admit an asymptotic expansion as 
$|\lambda| \to \infty$ and are obtained by a vector-valued  WKB ansatz, chosen in such a way that the error terms can be estimated in 
the most convenient way; see Section \ref{sec:special_solutions}, where we also prove the so called vanishing lemma. 
In Section \ref{sec:main_proof} we prove the above stated asymptotic estimates  as well as additional asymptotic estimates
for the norming constants $\kappa_n, \,n\in \mathbb{Z},$ introduced and studied in \cite{GK}, Section 8 and 10. 
The above stated results on the asymptotics of $\tau_n, \mu_n, \delta(\mu_n)$ and $\gamma_n$ are key ingredients in subsequent work to 
prove that the nonlinear Fourier transform of the defocusing NLS equation is semilinear \cite{KST2} and that the nonlinear part of
the solutions of the defocusing NLS equation on the circle is 1-smoothing \cite{KST3}.

\vspace{0.3cm}

\noindent{\em Related work:}
This paper is closely related to \cite{KST1} where asymptotic estimates of spectral quantities of Schr\"odinger operators $-\partial_x^2+q$ are
presented. In comparison with \cite{KST1}, notable differences are Theorem \ref{Theorem1.5} which will be used as a key ingredient in 
\cite{KST2}, as well as a conceptually new proof  of the asymptotic estimates of $\tau_n$ of Theorem \ref{Theorem1.4} (ii): 
Note that they cannot be obtained from Theorem \ref{Theorem1.2} (i); instead we derive them using Theorem \ref{Theorem1.4} (i) and 
Corollary \ref{Corollary1.3} together with the a priori estimate $\tau_n-\dot\lambda_n= O\left( \gamma_n^2\right),$ established in \cite{GK}, 
Lemma 6.9. The expansion of the eigenvalues of Sturm Liouville operators was pioneered by Marchenko \cite{Ma}. 
For selfadjoint ZS operators, the asymptotic estimates of the periodic eigenvalues of Theorem \ref{Theorem1.2} (ii) are stated 
(but not proved) in \cite{Ma}, p 94, except for the statement on the uniform boundedness of the error terms. 
Rough asymptotic estimates of $\gamma_n^2$ related to the problem of  characterizing the smoothness of $\varphi$ in terms of the decay of 
the $\gamma_n$ as $|n|\to \infty$ can be found in \cite{DM}, \cite{KaSeTo}, \cite{Ko} and \cite{Ma} as well as in references therein.
 
\vspace{0.3cm}

\noindent{\em Notation:} Throughout the paper we use for any $\lambda\in\mathbb{C}$, $x\in\mathbb{R}$,
and $\varphi\in L^2_c$ the following notation
\[
E_\lambda(x):= \begin{pmatrix}
e^{-i\lambda x}& 0\\ 0 &e^{i\lambda x}
\end{pmatrix},\quad R:= \begin{pmatrix}
i &0 \\ 0& -i
\end{pmatrix},\quad \Phi:= \begin{pmatrix}
0 & \varphi_1 \\ \varphi_2 & 0
\end{pmatrix} \,.
\]

\section{Special solutions}\label{sec:special_solutions}
In this section we prove estimates of special solutions $F_N\equiv F_N(x,\lambda) $ and $G_N\equiv G_N(x,\lambda)$ of 
the linear system $L(\varphi)F= \lambda F$, $\lambda\in \mathbb{C}\setminus\{0\},$ for $\varphi\in H^N_c$ with 
$N\in\mathbb{Z}_{\geq 1}$ which will be used to derive the asymptotics stated in the Introduction. 
These solutions are obtained with the WKB-type ansatz of the following form 
\begin{align}\label{3.1}
F_N(x,\lambda) := v_N(x,\lambda) 
\begin{pmatrix}
1 \\ \alpha_N(x,\lambda) 
\end{pmatrix} + \frac{R_N(x,\lambda)}{(2i\lambda)^N}
\end{align} while the error term $R_N(x,\lambda)$ satisfies 
$R_N(0,\lambda)= \begin{pmatrix}0\\0\end{pmatrix}$
and $v_N(x,\lambda)$ is the complex valued function 
\[v_N(x,\lambda) := \exp\Big(-i\lambda x+ i \int_0^x \varphi_1(t)\alpha_N(t,\lambda)\,dt\Big),\quad 
i\alpha_N(x,\lambda):= \sum_{n=1}^N\frac{r_n(x)}{(2i\lambda)^n}\]
and respectively, 
\begin{align}\label{3.2}
G_N(x,\lambda):= w_N(x,\lambda) \begin{pmatrix}
\beta_N(x,\lambda) \\1
\end{pmatrix} + \frac{S_N(x,\lambda)}{(2i\lambda)^N}
\end{align} where $S_N(0,\lambda)= \begin{pmatrix}0\\0\end{pmatrix}$,
\[
w_N(x,\lambda) := \exp\Big(i\lambda x- i \int_0^x \varphi_2(t)\beta_N(t,\lambda)\,dt\Big),\quad 
i\beta_N(x,\lambda):= \sum_{n=1}^N\frac{s_n(x)}{(2i\lambda)^n}.
\]
Substituting the ansatz \eqref{3.1} into $LF= \lambda F$ and using that 
$v_N'=(-i\lambda +i \varphi_1\alpha_N)v_N$ one gets 
\begin{align} \label{3.3}
(L-\lambda)\frac{R_N}{(2i\lambda)^N}= \begin{pmatrix}
0\\ \rho_N
\end{pmatrix}v_N
\end{align}
where
$\rho_N:= (i\alpha_N)'- 2i\lambda(i\alpha_N) +\varphi_1 (i\alpha_N)^2 -\varphi_2$ 
and $L\equiv L(\varphi)$.
The aim is to choose the coefficients $r_n(x)\equiv r_n(x,\varphi) $ in 
$i\alpha_N(x,\lambda) = \sum_{n=1}^N\frac{r_n(x)}{(2i\lambda)^n}$ so that all terms of $\rho_N(x,\lambda)$
of order $\le N-1$ in $1/\lambda$ vanish. 
We have,
\begin{align*}
\rho_N=& \sum_{n=1}^N\frac{r'_n(x)}{(2i\lambda)^n}-\sum_{n=0}^{N-1}\frac{r_{n+1}(x)}{(2i\lambda)^n}\\
& + \varphi_1\sum_{n=1}^N\Big(\sum_{k+ l=n}r_kr_l\Big)\frac{1}{(2i\lambda)^n}+ 
\varphi_1\sum_{n=N+1}^{2N}\Big(\sum_{k+l= n}r_kr_l\Big)\frac{1}{(2i\lambda)^n}-\varphi_2
\end{align*}
where we use the convention that the sums with lower limits greater than the upper limits vanish and that $r_n=0$ for
$n\le 0$ and $n\ge N+1$.
Collecting terms of the same order in $1/\lambda$ one gets  in the case $N=1$
\[
\rho_1=-(r_1+\varphi_2)+\frac{r_1'}{2i\lambda}+\frac{r_1^2\varphi_1}{(2i\lambda)^2}
\] and for $N\ge 2$
\begin{align*}
\rho_N=& -(r_1+\varphi_2)+(r_1'-r_2)\frac{1}{2i\lambda} + 
\sum_{n=2}^{N-1}\Big(r_n' -r_{n+1}+\varphi_1 \sum_{k=1}^{n-1}r_kr_{n-k}\Big)\frac{1}{(2i\lambda)^n}
\\ &+\Big( r_N'+ \varphi_1 \sum_{k=1}^{N-1}r_k r_{N-k}\Big)\frac{1}{(2i\lambda)^N} +
 \varphi_1\Big(\sum_{k=1}^{N}r_k r_{N+1-k}\Big)\frac{1}{(2i\lambda)^{N+1}} 
 \\& +\varphi_1\sum_{n=N+2}^{2N}\Big(\sum_{k=n-N}^{N}r_kr_{n-k}\Big)\frac{1}{(2i\lambda)^n}\,.
\end{align*}
For $\varphi\in H^N_c$ and $1\le n\le N$ we thus choose $r_1:=-\varphi_2,$ $r_2:= r_1'= -\varphi_2'$, and
\begin{align}
\label{3.3bis}
r_{n+1} := r_n'+\varphi_1 \sum_{k=1}^{n-1}r_kr_{n-k}\quad \forall 2\leq n\leq N-1\,.
\end{align}
This implies that,
\[
r_n=-\varphi_2^{(n-1)}+ p_n\,,
\] 
where $\varphi_2^{(n-1)}=\partial_x^{n-1}\varphi_2$, $p_1=p_2=0$, and for $3\leq n\leq N$,
$p_n$ is a polynomial in $\varphi_1$, $\varphi_2$ and its derivatives up to order $n-3$. 
Hence for any $1\leq n\leq N$,
\begin{align*}
r_n\in H^{N-n+1},\quad p_n\in H^{N-n+3}
\end{align*}
implying that
\begin{equation}\label{3.4'}
r_n\in H^1\,.
\end{equation}
Hence  $\alpha_N(x,\lambda)$ is a continuous function in $x$.
With this choice of $r_n$, $1\leq n\leq N$, $\rho_N$ can be written in the form
\begin{align*}
\rho_N= \frac{r_{N+1}}{(2i\lambda)^N}+ \varphi_1\Big(\sum_{k=1}^{N}r_k r_{N+1-k}\Big)\frac{1}{(2i\lambda)^{N+1}}
 +\varphi_1\!\!\!\sum_{n=N+2}^{2N}\!\!\!\Big(\sum_{k=n-N}^{N}r_k r_{n-k}\Big)\frac{1}{(2i\lambda)^n}
\end{align*}
where 
\begin{equation}\label{eq:r_{N+1}}
r_{N+1}:=r_N'+\varphi_1\sum_{k=1}^{N-1}r_k r_{N-k\,.}
\end{equation}
As above one sees that 
\[
r_{N+1}=-\varphi_2^{(N)}+p_{N+1}
\]
and $p_{N+1}\equiv p_{N+1}(\varphi)$ is a polynomial of $\varphi_1$, $\varphi_2$, and their derivatives up to order $N-2$. 
Hence,
\begin{align}\label{3.4bis}
r_{N+1}\in L^2\quad \text{and} \quad p_{N+1}\in H^2 \, .
\end{align}
By \eqref{3.4'} one has $\varphi_1\sum_{k=n-N}^{N}r_k r_{n-k}\in H^1$
for any $N+1\leq n\leq 2N$.
Hence $\rho_N$ is of the form
\begin{align}\label{3.5}
\rho_N= r_{N+1}\frac{1}{(2i\lambda)^N}+\tilde a_{N1}\frac{1}{(2i\lambda)^{N+1}}+ \tilde a_{N2} \frac{1}{(2i\lambda)^{N+2}} 
\end{align}
where ${\tilde a}_{N1}\equiv{\tilde a}_{N1}(\varphi):= \varphi_1 \sum_{k=1}^{N}r_k r_{N+1-k}\in H^1$ 
and ${\tilde a}_{N2}\equiv{\tilde a}_{N2}(\lambda, \varphi)$ is a polynomial in $1/\lambda$ of order $\le N-2$ with coefficients in $H^1$.
Equation \eqref{3.3} then reads 
\begin{align}\label{3.6}
(L-\lambda)R_N= \begin{pmatrix}
0\\ r_{N+1}+ {\tilde a}_{N1}\frac{1}{2i\lambda} + {\tilde a}_{N2}\frac{1}{(2i\lambda)^2}
\end{pmatrix}v_N
\end{align}
and 
\[
v_N(x,\lambda)=e^{-i x\lambda}\exp\left(\sum_{n=1}^N\Big(\int_0^x\varphi_1r_n\,dt\Big)
\frac{1}{(2i\lambda)^n}\right) \, .
\] 
By \eqref{3.4'}, $\int_0^x\varphi_1 r_n\,dt$ and consequently $v_N(\cdot,\lambda)$ are in $H^2([0,1],\mathbb{C})$. 
We have 
\begin{align}\label{3.7}
\exp\left(\sum_{n=1}^N\Big(\int_0^x\varphi_1r_n\,dt\Big)\frac{1}{(2i\lambda)^n}\right)=
1- \frac{1}{2i\lambda}\int_0^x \varphi_1\varphi_2\,dt+ O\Big(\frac{1}{\lambda^2}\Big)
\end{align}
where we use that $r_1= -\varphi_2$.
For any $\Lambda> 0,$ the estimate \eqref{3.7} holds uniformly for $|\lambda|\geq \Lambda,$ $0\leq x\leq 1,$ and 
uniformly on bounded sets of $\varphi$ in $H^N_c.$ Equation \eqref{3.6} then takes the form
\begin{align}\label{3.8}
(L-\lambda)R_N= E_{\lambda}(-x)\begin{pmatrix}
0\\ f_N
\end{pmatrix} ,\quad f_N= r_{N+1}+ a_{N1}\frac{1}{2i\lambda} + a_{N2}\frac{1}{(2i\lambda)^2}
\end{align}
where 
\[
a_{N1}(\varphi)= \varphi_1\sum_{k=1}^N r_k r_{N+1-k}-r_{N+1}\int_0^x \varphi_1\varphi_2\,dt \in L^2([0,1],\mathbb{C})
\] 
and $a_{N2}\equiv a_{N2}(\lambda,\varphi)$ is analytic as a map from $\mathbb{C}\setminus \{0\}\times H^N_c$ with values 
in $ L^2([0,1],\mathbb{C})$.
The maps $a_{N1}(\lambda,\varphi)$ and $a_{N2}(\lambda,\varphi)$ are bounded on bounded sets of $H^N_c$ uniformly in 
$|\lambda|\geq \Lambda.$

Now let us turn to the special solution $G_N.$
Substituting the ansatz \eqref{3.2} into the linear equation $LF= \lambda F$ and using that $w_N'= (i\lambda- i\varphi_2 \beta_N)w_N$ 
one gets 
\begin{align}\label{3.9}
(L-\lambda)\frac{S_N}{(2i\lambda)^N} = \begin{pmatrix}
\sigma_N \\ 0
\end{pmatrix}w_N
\end{align}
where $\sigma_N:= -(i\beta_N)'-2i\lambda(i\beta_N)+\varphi_2(i\beta_N)^2 -\varphi_1$.
Note that 
\begin{align*}
\sigma_N=& -\sum_{n=1}^N\frac{s'_n(x)}{(2i\lambda)^n}- \sum_{n=0}^{N-1}\frac{s_{n+1}(x)}{(2i\lambda)^n} \\
& + \varphi_2\sum_{n=2}^N\Big(\sum_{k+ l=n}s_ks_l\Big)\frac{1}{(2i\lambda)^n}  + 
\varphi_2 \sum_{n=N+1}^{2N}\Big(\sum_{k+l= n}s_ks_l\Big) \frac{1}{(2i\lambda)^n}-\varphi_1.
\end{align*}
Collecting terms of the same order in $1/\lambda$ one gets in the case $N=1$
\[
\sigma_1=-(s_1+\varphi_1)-\frac{s_1'}{2 i\lambda}+\frac{s_1^2\varphi_2}{(2 i\lambda)^2},
\]
and for $N\ge 2$,
\begin{align*}
\sigma_N=& -(s_1+\varphi_1)-(s_1'+s_2)\frac{1}{2i\lambda}
+\sum_{n=2}^{N-1}\Big(-s_n' -s_{n+1}+\varphi_2 \sum_{k=1}^{n-1}s_ks_{n-k}\Big)\frac{1}{(2i\lambda)^n}\\ &
+\Big( -s_N'+ \varphi_2 \sum_{k=1}^{N-1}s_ks_{N-k}\Big)\frac{1}{(2i\lambda)^N}
+\varphi_2\Big(\sum_{k=1}^{N}s_ks_{N+1-k}\Big)\frac{1}{(2i\lambda)^{N+1}} \\& 
+\varphi_2\sum_{n=N+2}^{2N}\Big(\sum_{k=n-N}^{N}s_ks_{n-k}\Big)\frac{1}{(2i\lambda)^n}.
\end{align*}
For $\varphi\in H^N_c$ and $1\le n\le N$ we choose $s_1:=-\varphi_1,$ $s_2:= -s_1'= \varphi_1'$, and
\begin{align}\label{3.9bis}
s_{n+1} := -s_n'+\varphi_2 \sum_{k=1}^{n-1}s_ks_{n-k}\qquad \forall 2\leq n\leq N-1\,.
\end{align}
This implies that,
\[
s_n=(-1)^n\varphi_1^{(n-1)}+ q_n,
\]
where $\varphi_1^{(n-1)}=\partial_x^{n-1}\varphi_1$, $q_1=q_2=0$, and for $3\leq n\leq N$, $q_n$ is a polynomial in 
$\varphi_1$, $\varphi_2$ and its derivatives up to order $n-3$. 
Hence, for any $1\leq n\leq N$,
\[
s_n\in H^{N-n+1}\quad \text{and} \quad q_n\in H^{N-n+3}
\]
implying that
\begin{equation}\label{3.10}
s_n\in H^1\,.
\end{equation}
Hence, $\beta_N(x,\lambda)$ is a continuous function of $x$.
With this choice of $s_n,\, 1\leq n\leq N$, $\sigma_N$ can be written in the form
\begin{align*}
\sigma_N = \frac{s_{N+1}}{(2i\lambda)^N}+ \varphi_2 \Big(\sum_{k=1}^{N}s_ks_{N+1-k}\Big)\frac{1}{(2i\lambda)^{N+1}}
 +\varphi_2\!\!\!\!\!\sum_{n=N+2}^{2N}\!\!\!\!\Big(\sum_{k=n-N}^{N}s_k s_{n-k}\Big)\frac{1}{(2 i\lambda)^n}
\end{align*}
where
\begin{equation}\label{eq:s_{N+1}}
s_{N+1}:=-s_N'+\varphi_2\sum_{k=1}^{N-1}s_k s_{N-k}\,.
\end{equation}
As above one sees that 
\[
s_{N+1}=(-1)^{N+1}\varphi_1^{(N)}+q_{N+1}
\]
where $q_{N+1}\equiv q_{N+1}(\varphi)$ is a polynomial of $\varphi_1$, $\varphi_2$, and their derivatives up to order $N-2$.
Hence, 
\begin{align}\label{3.10bis}
s_{N+1}\in L^2\quad \text{and} \quad q_{N+1}\in H^2.
\end{align}
By \eqref{3.10} one has $\varphi_2\sum_{k=n-N}^{N}s_ks_{n-k}\in H^1$
for any $N+1\leq n\leq 2N$. 
Hence $\sigma_N$ is of the form
\begin{align}\label{3.11}
\sigma_N= s_{N+1}\frac{1}{(2i\lambda)^N}+{\tilde b}_{N1}\frac{1}{(2i\lambda)^{N+1}}+ {\tilde b}_{N2} \frac{1}{(2i\lambda)^{N+2}} 
\end{align}
where $\tilde b_{N1}\equiv{\tilde b}_{N1}(\varphi):= \varphi_2 \sum_{k=1}^{N}s_ks_{N+1-k}\in H^1$ and 
${\tilde b}_{N2}\equiv \tilde b_{N2}(\lambda, \varphi)$ is a polynomial in $1/\lambda$ of order $\le N-2$ with coefficients in $H^1$.
Equation \eqref{3.9} then reads 
\begin{align}\label{3.12}
(L-\lambda)S_N= \begin{pmatrix}
 s_{N+1}+ \tilde b_{N1}\frac{1}{2i\lambda} + \tilde b_{N2}\frac{1}{(2i\lambda)^2}\\0
\end{pmatrix}w_N
\end{align}
with 
\[
w_N(x,\lambda)=e^{i x\lambda}\exp\left(-\sum_{n=1}^N\Big(\int_0^x\varphi_2s_n\,dt\Big)\frac{1}{(2i\lambda)^n}\right)\,.
\]
By \eqref{3.10}, $\int_0^x\varphi_2 s_n\,dt$ and consequently $w_N(\cdot,\lambda)$ are in $H^2([0,1],\mathbb{C})$.
Furthermore,
\begin{align}\label{3.13}
\exp\left(-\sum_{n=1}^N\Big(\int_0^x\varphi_2s_n\,dt\Big)\frac{1}{(2i\lambda)^n}\right)
=1+ \frac{1}{2 i\lambda}\int_0^x \varphi_1\varphi_2\,dt+ O\Big(\frac{1}{\lambda^2}\Big)
\end{align}
where we use that $s_1= -\varphi_1$.
For any $\Lambda> 0,$ the estimate \eqref{3.13} holds uniformly for $|\lambda|\geq \Lambda,$ $0\leq x\leq 1$,
and uniformly on bounded sets of $\varphi$ in $H^N_c$.
Then \eqref{3.12} reads
\begin{align}\label{3.13bis}
(L-\lambda)S_N= E_{\lambda}(-x)\begin{pmatrix}
g_N\\ 0
\end{pmatrix} ,\quad g_N= s_{N+1}+ b_{N1}\frac{1}{2i\lambda} + b_{N2}\frac{1}{(2i\lambda)^2},
\end{align}
where 
\[
b_{N1}(\varphi) = \varphi_2\sum_{k=1}^N s_ks_{N+1-k}+ s_{N+1}\int_0^x \varphi_1\varphi_2\,dt \in L^2([0,1],\mathbb{C}) 
\] 
and $ b_{N2}(\lambda,\varphi)$ is analytic as a map from $\mathbb{C}\setminus \{0\}\times H^N_c$ with values in $L^2([0,1],\mathbb{C}) $.
The maps $b_{N1}(\varphi)$ and $b_{N2}$ are bounded on bounded sets of $H^N_c$ uniformly in $|\lambda|\geq \Lambda.$

As a next step we want to estimate $R_N\,[S_N]$ by using that it satisfies the inhomogeneous linear ODE \eqref{3.8} [\eqref{3.13bis}] with
initial conditions 
$R_N(0,\lambda)=\begin{pmatrix}
0\\0
\end{pmatrix}\,\big[S_N(0,\lambda)=\begin{pmatrix}
0\\0
\end{pmatrix}\big].$ 

For the proof of the main results in Section \ref{sec:main_proof} we need the asymptotics  of $R_N(1,\xi_n)$
and $S_N(1,\xi_n)$ as $|n|\to \infty$ for sequences $\xi_n= n\pi +O\big(\frac{1}{n}\big) .$
Denote $\langle n\rangle:=1+|n|$.

\begin{prop1}\label{Proposition3.2}
For a given sequence $(\xi_n)_{n\in\mathbb{Z}}$ of complex numbers $\xi_n=n\pi +\alpha_n$ such that
$|\alpha_n|\leq \frac{a}{\langle n\rangle}$ for some positive (independent of $n$) constant $a>0$
and for any $\varphi\in H^N_c$
\begin{itemize}
\item[(i)]
$ R_N(1,\xi_n)= \begin{pmatrix}
\frac{(-1)^n}{2i\xi_n}\int_0^1\varphi_1r_{N+1}\,dt \\
i(-1)^{n+1} \widehat{(\varphi_2^{(N)})}(n)
\end{pmatrix}+ \frac{\ell^2_n}{n}
$
\item[(ii)]
$  S_N(1,\xi_n)= \begin{pmatrix}
i (-1)^{n}(-1)^N \widehat{(\varphi_1^{(N)})}(-n)\\
\frac{(-1)^{n+1}}{2i\xi_n}\int_0^1\varphi_2 s_{N+1}\,dt 
\end{pmatrix}+
\frac{\ell^2_n}{n}
$
\end{itemize}
where the estimates hold uniformly for $(\xi_n)_{n\in\mathbb{Z}}$ with 
$|\alpha_n|\leq \frac{a}{\langle n\rangle}$ and uniformly on bounded sets of $\varphi$'s in $H^N_c.$
\end{prop1}
\begin{proof}
The estimates (i) and (ii) are proved in a similar way and so we concentrate on (i) only.
Recall that by \eqref{3.8}, $(L-\xi_n)R_N= E_{\xi_n}(-x)
\begin{pmatrix}
0\\ f_N
\end{pmatrix}$
where $f_N= r_{N+1}+a_{N1}\frac{1}{2i\xi_n}+a_{N2}\frac{1}{(2i\xi_n)^2}.$
In addition, $R_N$ satisfies $R_N(0,\xi_n)= (0,0).$ 
As $e^{i \xi_n s}=O(1)$ uniformly in $s\in[-1,1]$ and $n\in\mathbb{Z}$, we get from
Corollary \ref{CorollaryA.3} in Appendix A, with $M_1(x,\lambda)$ and $Q(x)$ defined as in \eqref{A.1} and \eqref{A.2}
respectively, that $R_N(1,\xi_n)$ admits the following asymptotic expansion as $|n|\to \infty$ 
\[
R_N(1,\xi_n)  = \sum_{k=1}^2 A_{Nk}^R(1,\xi_n)+\frac{1}{2\xi_n}\sum_{k=1}^4 B_{Nk}^R(1,\xi_n)+ O\Big(\frac{1}{n^2}\Big)
\]
where 
\begin{align*}
A_{N1}^R(1,\xi_n)=& \begin{pmatrix}
0 \\ i\int_0^1e^{i\xi_n(1-2t)}r_{N+1}(t)\,dt
\end{pmatrix}
\end{align*}
\begin{align*}
A_{N2}^R(1,\xi_n)=& \begin{pmatrix}
0 \\ \frac{1}{2\xi_n}\int_0^1e^{i\xi_n(1-2t)}a_{N1}(t)\,dt
\end{pmatrix}
\end{align*}
\begin{align*}
B_{N1}^R(1,\xi_n)=&2i\xi_n M_1(1,\xi_n)\begin{pmatrix}
0\\
\int_0^1 e^{-2i\xi_nt}r_{N+1}(t)\,dt
\end{pmatrix}
\end{align*}
\begin{align*}
B_{N2}^R(1,\xi_n)=&\begin{pmatrix}
0\\
Q(1)\int_0^1 e^{i\xi_n(1-2t)}r_{N+1}(t)\,dt
\end{pmatrix}
\end{align*}
\begin{align*}
B_{N3}^R(1,\xi_n)=&-iE_{\xi_n}(1)\int_0^12\xi_n M_1(t,\xi_n)\begin{pmatrix}
0\\
 e^{-i\xi_n t}r_{N+1}(t)
\end{pmatrix}\,dt
\end{align*}
\begin{align*}
\\ B_{N4}^R(1,\xi_n)=&
\begin{pmatrix}
0 \\ -\int_0^1e^{i\xi_n(1-2t)}Q(t)r_{N+1}(t)\,dt
\end{pmatrix}
\end{align*}
and the estimate is uniform on bounded sets of $H^N_c$ and on sequences $\xi_n= n\pi +\alpha_n,n\in \mathbb{Z},$ 
with $|\alpha_n|\leq \frac{a}{\langle n\rangle}.$
The terms in the expansion are treated individually. Concerning $A^R_{N1}(1,\xi_n),$ recall that $r_{N+1}=-\varphi_2^{(N)}+p_{N+1}$ 
where $p_{N+1}\in H^2.$ By Lemma \ref{LemmaB.1} of Appendix B, 
\[
\int_0^1 e^{i\xi_n(1-2t)}\varphi_2^{(N)}(t)\,dt 
= (-1)^n\widehat{(\varphi_2^{(N)})}(n)+\frac{\ell^2_n}{n}.
\]
Furthermore, integrating by parts and using that $p_{N+1}$ is $1$--periodic 
\begin{align*}
\int_0^1 e^{i\xi_n(1-2t)}p_{N+1}(t)\,dt = & \frac{1}{2i\xi_n}p_{N+1}(0) 2i \sin\xi_n  + 
\frac{1}{2i\xi_n}\int_0^1 e^{i\xi_n(1-2t)}p_{N+1}'(t)dt \, .
\end{align*}
Note that $\sin\xi_n= (-1)^n\sin \alpha_n= O\big(\frac{1}{n}\big).$ Integrating by parts once again then yields 
\[
\int_0^1 e^{i\xi_n(1-2t)}p_{N+1}(t)dt =O\Big(\frac{1}{n^2}\Big).\]
Altogether we thus have proved that
\begin{align*}A_{N1}^R(1,\xi_n)=& \begin{pmatrix}
0 \\ i(-1)^{n+1}\widehat{(\varphi_2^{(N)})}(n) + \frac{\ell^2_n }{n}
\end{pmatrix} \, .
\end{align*}
Towards $A_{N2}^R(1,\xi_n),$ recall that by \eqref{3.8}, 
\[a_{N1}= \varphi_1\sum_{k=1}^{N}  r_kr_{n+1-k}-r_{N+1}Q\in L^2([0,1],\mathbb{C}) \, .
\]
By Lemma \ref{LemmaB.1} in Appendix B it then follows that 
\begin{align*}
A_{N2}^R(1,\xi_n)=& \begin{pmatrix}
0 \\ \frac{1}{2\xi_n}\int_0^1e^{i\xi_n(1-2t)}a_{N1}(t)dt
\end{pmatrix}= \begin{pmatrix}
0\\\frac{\ell^2_n}{n}
\end{pmatrix} \, .
\end{align*}
Concerning 
$B_{N1}^R(1,\xi_n),$ recall that by \eqref{A.1} in Appendix A,
\[
2\xi_n M_1(1,\xi_n) = \big(E_{\xi_n}(-1)-E_{\xi_n}(1)\big)\Phi(0)- E_{\xi_n}(1)P_{\xi_n}(1)
\]
where $P_{\xi_n}(1)= \int_0^1E_{\xi_n}(-2t) \Phi'(t)\,dt$.
Hence $2\xi_n M_1(1,\xi_n) =O(1)$ and using that $M_1(1,\xi_n)$ is off-diagonal one concludes again from 
Lemma \ref{LemmaB.1} in Appendix B that
\begin{align*}
B_{N1}^R(1,\xi_n)=&2i\xi_n M_1(1,\xi_n)\begin{pmatrix}
0\\
\int_0^1 e^{-2i\xi_nt}r_{N+1}(t)dt
\end{pmatrix}=\begin{pmatrix}
\ell^2_n\\ 0
\end{pmatrix}.
\end{align*}
Similarly, one sees that 
\begin{align*}
B_{N2}^R(1,\xi_n)=&\begin{pmatrix}
0\\
Q(1)\int_0^1 e^{i\xi_n(1-2t)}r_{N+1}(t)dt
\end{pmatrix} = \begin{pmatrix}0\\
\ell^2_n
\end{pmatrix}
\end{align*}
and
\begin{align*} B_{N4}^R(1,\xi_n)=&
\begin{pmatrix}
0 \\ -\int_0^1e^{i\xi_n(1-2t)}Q(t)r_{N+1}(t)\,dt\end{pmatrix}=\begin{pmatrix}0\\
\ell^2_n 
\end{pmatrix} \, .
\end{align*}
It remains to consider 
\begin{align*}
B_{N3}^R(1,\xi_n)=&-i E_{\xi_n}(1)\int_0^12\xi_n M_1(t,\xi_n)\begin{pmatrix}
0\\
 e^{-i\xi_n t}r_{N+1}(t)
\end{pmatrix}dt \, .
\end{align*}
By \eqref{A.1} 
\[
2\xi_n M_1(t,\xi_n)= E_{\xi_n}(-t)\Phi(t)- E_{\xi_n}(t)\Phi(0)- E_{\xi_n}(t) P_{\xi_n}(t)
\]
where $P_{\xi_n}(t)= \int_0^t E_{\xi_n}(-2x)\Phi'(x)dx.$ 
Hence
\begin{align*}&2\xi_n M_1(t,\xi_n)\begin{pmatrix}
0\\ e^{-i\xi_n t}r_{N+1}(t)
\end{pmatrix}\\& =
\begin{pmatrix}
\varphi_1(t)r_{N+1}(t)-e^{2i\xi_n t} \varphi_1(0) r_{N+1}(t) - r_{N+1}(t)\int_0^t e^{i 2\xi_n(x-t)}\varphi_1'(x)\,dx\\0
\end{pmatrix}.
\end{align*}
By Lemma \ref{LemmaB.1} in Appendix B, $\int_0^1 e^{-2i\xi_n t}r_{N+1}(t)\in \ell^2_n$ and 
\[
\int_0^1 e^{-2i\xi_n t}\Big(r_{N+1}(t)\int_0^t e^{i2\xi_n x}\varphi_1'(x)\,dx\Big)\,dt\in \ell^2_n.
\]
As $e^{-i\xi_n}=(-1)^n+ O\big(\frac{1}{n}\big)$ it then follows that 
\[
B_{N3}^R(1,\xi_n)=\begin{pmatrix}-i(-1)^n \int_0^1 \varphi_1(t)r_{N+1}(t) dt\\ 0\end{pmatrix} +
\begin{pmatrix}\ell^2_n\\0\end{pmatrix}.
\]
Altogether we have proved that the claimed asymptotics of $R_N(1,\xi_n)$. 
Going through the arguments of the proof one verifies that the claimed uniformity statement holds. 
\end{proof}

\noindent Next we will prove the following vanishing lemma. 
\begin{lem1}\label{Lemma3.1}
For any $\varphi \in H^N_c$ with $N\in \mathbb{Z}_{\geq 1},$ 
\[\int_0^1\varphi_1 r_k\,dt=\int_0^1\varphi_2s_k\,dt\qquad \forall 1\leq k\leq N+1\]
where $r_k\,[s_k]$ are given by \eqref{3.3bis} and \eqref{eq:r_{N+1}} 
\mbox{\rm{[}}\eqref{3.9bis} and \eqref{eq:s_{N+1}}\mbox{\rm{]}}.
As a consequence, $\int_0^1 \varphi_1(t)\alpha_N(t,\lambda)\,dt=\int_0^1 \varphi_2(t)\beta_N(t,\lambda)\,dt $ for any 
$\lambda\in \mathbb{C}.$
\end{lem1}
\begin{proof}
Let $\varphi$ be an arbitrary element of $H^N_c$ with $N\in \mathbb{Z}_{\geq 1}$ and let $\xi_n$, $n\in\mathbb{Z}$, be as
in Proposition \ref{Proposition3.2}. 
The claimed identities follow from the Wronskian identity, applied to the special solutions $F_N,\,G_N$ constructed above, 
\begin{align}\label{3.14}
\det \big[F_N(1,\xi_n)\,\,G_N(1,\xi_n)\big]=\det \big[F_N(0,\xi_n)\,\,G_N(0,\xi_n)\big],\,\,\,n\in\mathbb{Z}.
\end{align}  
By the definition of $F_N$ and $G_N$ one has 
\[
\det \big[F_N(0,\xi_n)\,\,G_N(0,\xi_n)\big]= 
\det\begin{pmatrix}
1& \beta_N(0,\xi_n) \\
\alpha_N(0,\xi_n) &1
\end{pmatrix}.
\]
To compute the left hand side of \eqref{3.14}, note that by \eqref{3.1} and \eqref{3.2}, 
\[
F_N(1,\xi_n)= e^{-i\xi_n} 
\begin{pmatrix}
1\\ \alpha_N(0,\xi_n)
\end{pmatrix} 
\exp\Big(i \int_0^1\varphi_1(t)\alpha_N(t,\xi_n),dt\Big)+\frac{ R_N(1,\xi_n)}{(2i\xi_n)^N}
\]
and 
\[
G_N(1,\xi_n)= e^{i\xi_n} 
\begin{pmatrix}
\beta_N(0,\xi_n)\\1
\end{pmatrix} 
\exp\Big(-i \int_0^1\varphi_2(t)\beta_N(t,\xi_n)\,dt\Big)+\frac{S_N(1,\xi_n)}{(2i\xi_n)^N}
\] 
where we used that $\alpha_N(x,\xi_n)$ and $\beta_N(x,\xi_n)$ are both $1$--periodic in $x$. 
This together with Proposition \ref{Proposition3.2} imply that $\det \big[F_N(1,\xi_n)\,\,G_N(1,\xi_n)\big]$
satisfies the estimate
\[  
e^{i \int_0^1\big(\varphi_1(t)\alpha_N(t,\xi_n)-\varphi_2(t)\beta_N(t,\xi_n)\big)\,dt}
\det\begin{pmatrix}
1& \beta_N(0,\xi_n) \\
\alpha_N(0,\xi_n) &1
\end{pmatrix} + O\Big(\frac{1}{n^{N+1}}\Big)\,.
\]
For $|n|$ sufficiently large, $|\alpha_N(0,\xi_n)\beta_N(0,\xi_n)|\leq \frac{1}{2}$ and
hence 
\[
\Big|
\det\begin{pmatrix}
1& \beta_N(0,\xi_n) \\
\alpha_N(0,\xi_n) &1
\end{pmatrix}
\Big|\ne 0 \, ,
\]
implying that $\exp\Big(i \int_0^1\big(\varphi_1\alpha_N(t,\xi_n)-\varphi_2\beta_N(t,\xi_n)\big)\,dt\Big)= 
1+O\big(\frac{1}{n^{N+1}}\big)$
or 
\[
\exp\Big(\sum_{k=1}^N\frac{1}{(2i\xi_n)^k} \int_0^1\big(\varphi_1(t)r_k(t)-\varphi_2(t)s_k(t)\big)\,dt\Big)= 
1+O\Big(\frac{1}{n^{N+1}}\Big).
\]
Taking the logarithm of both sides of this formula for $n$ sufficiently large one concludes that 
\begin{align}\label{3.15}
\int_0^1\varphi_1(t)r_k(t)dt=\int_0^1\varphi_1(t)s_k(t)\,dt\quad \forall
1\leq k\leq N. 
\end{align}
In case $\varphi\in H^{N+1}_c,$ the latter identity also holds for $k=N+1$. As $r_{N+1}$ and $s_{N+1}$ are polynomials in
$\varphi_1$, $\varphi_2$ and their derivatives up to order $N$, the identity continues to hold for any $\varphi \in H^N_c$ as 
the embedding $H^{N+1}_c\hookrightarrow H^N_c$ is dense. 
\end{proof}
 
\noindent Lemma \ref{Lemma3.1} and Proposition \ref{Proposition3.2} lead to the following formulas for $F_N(1,\lambda)$ and 
$G_N(1,\lambda),$\begin{align}
\label{3.16}
F_N(1,\lambda)= &\begin{pmatrix}
1\\\alpha_N(0,\lambda)
\end{pmatrix} e^{-i\theta_N}+\frac{1}{(2i\lambda)^N} R_N(1,\lambda)\\
\label{3.17}
G_N(1,\lambda)=& \begin{pmatrix}
\beta_N(0,\lambda)\\ 1
\end{pmatrix} e^{i\theta_N}+\frac{1}{(2i\lambda)^N} S_N(1,\lambda)
\end{align}
where 
\[
\alpha_N(0,\lambda)=- i \sum_{k=1}^N\frac{r_k(0)}{(2 i\lambda)^k},\,\,\,\,
\beta_N(0,\lambda)=- i \sum_{k=1}^N\frac{s_k(0)}{(2 i\lambda)^k}, 
\]
\begin{align}\label{3.17bis} 
\theta_N(\lambda):= \lambda- \int_0^1\varphi_1\alpha_N dt= 
\lambda+i\sum_{k=1}^N\Big(\int_0^1\varphi_1(t)r_k(t)\,dt\Big)\frac{1}{(2i\lambda)^k}\,.
\end{align}
Furthermore, in view of Lemma \ref{Lemma3.1}, the estimate of $S_N(1,\xi_n)$ of Proposition \ref{Proposition3.2} can be
written in terms of $r_{N+1}$ instead of $s_{N+1}$.
The two estimates of Proposition \ref{Proposition3.2} thus read
\[
R_N(1,\xi_n)= \begin{pmatrix}
\frac{(-1)^n}{2i\xi_n}\int_0^1\varphi_1r_{N+1}\,dt \\
i(-1)^{n+1} \widehat{(\varphi_2^{(N)})}(n)
\end{pmatrix}+ \frac{\ell^2_n}{n}
\]
\[
S_N(1,\xi_n)= \begin{pmatrix}
i (-1)^{n}(-1)^N \widehat{(\varphi_1^{(N)})}(-n)\\
\frac{(-1)^{n+1}}{2i\xi_n}\int_0^1\varphi_1 r_{N+1}\,dt 
\end{pmatrix}+
\frac{\ell^2_n}{n}\,.
\]
Here we used that $\alpha_N(1,\lambda)=\alpha_N(0,\lambda)$ and $\beta_N(1,\lambda)=\beta_N(0,\lambda)$.

\medskip\medskip

\noindent As an application of the estimates obtained so far, we consider the $2\times 2$ matrix 
$M_N(x,\lambda):= 
\begin{pmatrix}
F_N(x,\lambda)\,\,G_N(x,\lambda)
\end{pmatrix}$
with columns $F_N(x,\lambda)$ and $G_N(x,\lambda).$ It follows from the definition of $\alpha_N$ and $\beta_N$ that 
there exists $\Lambda>0$ so that $|\alpha_N(0,\lambda)\beta(0,\lambda)|\leq 1/2$ implying that in view of the 
Wronskian identity, for any $|\lambda|\geq \Lambda$ and $0\leq x\leq 1$
\[
|\det M_N(x,\lambda)|= |\det M_N(0,\lambda)|= |1-\alpha_N(0,\lambda)\beta(0,\lambda)|\geq 1/2.
\]
The constant $\Lambda$ can be chosen uniformly on bounded sets of $\varphi$'s in $H^N_c.$ 
For $|\lambda|\geq \Lambda$ one then verifies that $M(x,\lambda)= M_N(x,\lambda) M_N(0,\lambda)^{-1}$ reads
\begin{align}\label{3.18}
M(x,\lambda)=  \frac{1}{1-\alpha_N^0\beta_N^0}
\begin{pmatrix}
F_N-\alpha_N^0G_N & G_N-\beta_N^0F_N 
\end{pmatrix}
\end{align}
where $\alpha_N^0\equiv \alpha_N^0(\lambda):=\alpha_N(0,\lambda)$
and $\beta_N^0\equiv \beta_N^0(\lambda):=\beta_N(0,\lambda).$
For $\xi_n= n\pi+\alpha_n$, with $|\alpha_n|\le\frac{a}{\langle n\rangle}$, $a>0$, and $|\xi_n|\geq \Lambda,$ one then gets
\begin{align}\label{3.18bis}
M(1,\xi_n)=\frac{1}{1-\alpha_N^0(\xi_n)\beta_N^0(\xi_n)}
\begin{pmatrix}  \grave{m}_{N1}(\xi_n) &\grave{m}_{N2}(\xi_n)\\
\grave{m}_{N3}(\xi_n) &\grave{m}_{N4}(\xi_n)
\end{pmatrix}
\end{align}
where 
\begin{align*}
\grave{m}_{N1}(\xi_n):=& e^{-i\theta_N(\xi_n)}-\alpha_N^{0}(\xi_n)\beta_N^{0}(\xi_n)e^{i\theta_N(\xi_n)} \\ &+
\frac{(-1)^n}{(2i\xi_n)^{N+1}}\int_0^1\varphi_1 r_{N+1}\,dt +\frac{\ell^2_n}{n^{N+1}}
\\
\grave{m}_{N2}(\xi_n):=& \beta_N^{0}(\xi_n) \big(e^{i\theta_N(\xi_n)}- e^{-i\theta_N(\xi_n)}\big) +
\frac{i(-1)^{n}}{(-2i\xi_n)^{N}}\widehat{(\varphi_1^{(N)})}(-n) +\frac{\ell^2_n}{n^{N+1}}
\\
\grave{m}_{N3}(\xi_n):=& \alpha_N^{0}(\xi_n)\big(e^{-i\theta_N(\xi_n)}-e^{i\theta_N(\xi_n)}\big) +
\frac{i (-1)^{n+1}}{(2i\xi_n)^{N}}\widehat{(\varphi_2^{(N)})}(n) +\frac{\ell^2_n}{n^{N+1}}
\\
\grave{m}_{N4}(\xi_n):=& e^{i\theta_N(\xi_n)}-\alpha_N^{0}(\xi_n)\beta_N^{0}(\xi_n)e^{-i\theta_N(\xi_n)}+\\&
+\frac{(-1)^{n+1}}{(2i\xi_n)^{N+1}}\int_0^1\varphi_1r_{N+1}\,dt +\frac{\ell^2_n}{n^{N+1}}
\end{align*}
Recall that $\Delta(\lambda)\,[\delta(\lambda)]$ denotes the trace [anti--trace] of the Floquet matrix $M(1,\lambda)$ 
whereas 
\[
2i\chi_D(\lambda)= (m_4+m_3-m_2-m_1)\vert_{(1,\lambda)}.
\] 
We obtain the following 
\begin{prop1}\label{Corollary3.3}
Let $\varphi\in H^N_c$ and let $(\xi_n)_{n\in\mathbb{Z}}$ be a sequence of complex numbers $\xi_n= n\pi+\alpha_n$ such that 
$|\alpha_n|\leq \frac{a}{\langle n\rangle}$ for some positive (independent of $n$) constant $a>0$.
Then for $|n|$ sufficiently large so that $|\alpha_N^{0}(\xi_n)\beta_N^{0}(\xi_n)|\leq 1/2$, the following holds:
\begin{itemize}
\item[(i)] $\Delta(\xi_n)= 2\cos \theta_N(\xi_n)+ \frac{\ell^2_n}{n^{N+1}}$;
\item[(ii)] 
$
\delta(\xi_n)=\frac{\beta_N^{0}(\xi_n)-\alpha_N^{0}(\xi_n)}{1-\alpha_N^{0}(\xi_n)\beta_N^{0}(\xi_n)}  2i\sin \theta_N(\xi_n) +
i (-1)^n\big(\hat\varphi_1(-n)-\hat\varphi_2(n)\big)+ \frac{\ell^2_n}{n^{N+1}};
$
\item[(iii)]
\begin{align*}
2i\chi_D(\xi_n)=& \frac{\big(1-\alpha_N^{0}(\xi_n)\big)\big(1-\beta_N^{0}(\xi_n)\big)}{1-\alpha_N^{0}(\xi_n)\beta_N^{0}(\xi_n)} 
2i\sin \theta_N(\xi_n)+\\&+ i (-1)^{n+1}\big(\hat\varphi_1(-n)+\hat\varphi_2(n)\big)+\\&+ 
\frac{2(-1)^{n+1}}{(2i\pi\xi_n)^{N+1}}\int_0^1\varphi_1r_{N+1}\,dt + \frac{\ell^2_n}{n^{N+1}}.
\end{align*} 
\end{itemize}
These estimates hold uniformly for $(\xi_n)_{n\in\mathbb{Z}}$ with $|\alpha_n|\le\frac{a}{\langle n\rangle}$ and uniformly
on bounded sets of $\varphi$'s in $H^N_c$.
\end{prop1}
We finish this section by providing asymptotic expansions for $\dot R_N(1,\xi_n)= \partial_\lambda R_N(1,\xi_n)$ and
$\dot S_N(1,\xi_n)= \partial_\lambda S_N(1,\xi_n).$ Arguing as in the proof of Proposition \ref{Proposition3.2} one obtains the following 
\begin{prop1}\label{Proposition3.4}
For complex numbers $\xi_n=n\pi +\alpha_n$ with $|\alpha_n|\leq \frac{a}{\langle n\rangle}$
and $\varphi\in H^N_c$
\begin{itemize}
\item[(i)] 
$ \dot R_N(1,\xi_n)= \begin{pmatrix}
\frac{i(-1)^{n+1}}{2i\xi_n}\int_0^1\varphi_1(t)r_{N+1}(t)\,dt \\
(-1)^{n}\int_0^1 e^{-i2n\pi t}(2t-1)r_{N+1}(t)\,dt
\end{pmatrix}+ \frac{\ell^2_n}{n}$
\item[(ii)] 
$ \dot S_N(1,\xi_n)= \begin{pmatrix}
(-1)^{n}\int_0^1 e^{i2n\pi t}(2t-1)s_{N+1}(t)\,dt\\
\frac{i(-1)^{n+1}}{2i\xi_n}\int_0^1\varphi_1(t)r_{N+1}(t)\,dt
\end{pmatrix}+ \frac{\ell^2_n}{n}$
\end{itemize}
where the estimates hold uniformly for $(\xi_n)_{n\in\mathbb{Z}}$ with $|\alpha_n|\leq \frac{a}{\langle n\rangle}$ and 
uniformly on bounded sets of $\varphi$'s in $H^N_c.$
\end{prop1} 
Proposition \ref{Proposition3.4} leads to the following asymptotics for $\dot\Delta(\xi_n)  .$
\begin{cor1}\label{Corollary3.5}
For complex numbers $\xi_n=n\pi +\alpha_n$ with $|\alpha_n|\leq \frac{a}{\langle n\rangle}$ and $\varphi\in H^N_c$
\begin{align*}
\dot\Delta(\xi_n)=& -\dot \theta_N(\xi_n)2 \sin \theta_N(\xi_n)+ \frac{2i(-1)^{n+1}}{(2i\xi_n)^{N+1}}\int_0^1
\varphi_1(t)r_{N+1}(t)\,dt+\frac{1}{n^{N+1}}\ell^2_n
\end{align*}
where the estimate holds uniformly for $(\xi_n)_{n\in\mathbb{Z}}$ with $|\alpha_n|\leq \frac{a}{\langle n\rangle}$ and 
and uniformly on bounded sets of $\varphi$'s in $H^N_c.$
\end{cor1}
\begin{proof}
In view of \eqref{3.18bis} we have 
\begin{align}
\nonumber
\Delta(\lambda)=& 2\cos \theta_N(\lambda)+ \frac{1}{1-\alpha_N^0\beta_N^0}\frac{1}{(2i\lambda)^N}
\Big(\big(R_N(1,\lambda)-\alpha_N^0S_N(1,\lambda)\big)_1\\\label{3.30}
& +\big (S_N(1,\lambda)-\beta_N^0 R_N(1,\lambda)\big)_2\Big)
\end{align}
where we denoted by $(\cdot)_1\,[(\cdot)_2]$ the first [second] component of 
the expression $R_N(1,\lambda)-\alpha_N^0S_N(1,\lambda)$ $[S_N(1,\lambda)-\beta_N^0R_N(1,\lambda)].$
Recall that 
\begin{align*}
\alpha_N^0=& \alpha_N(0,\lambda) = -i \sum_{k=1}^N\frac{1}{(2i\lambda)^k}r_k(0)= O\Big(\frac{1}{\lambda}\Big)
\\\beta_N^0=& \beta_N(0,\lambda) = -i \sum_{k=1}^N\frac{1}{(2i\lambda)^k}s_k(0)= O\Big(\frac{1}{\lambda}\Big)
\end{align*}
implying that $(\alpha_N^0)^{\cdot}, (\beta_N^0)^{\cdot}= O\big(\frac{1}{\lambda^2}\big)$ and  
$(\alpha_N^0\beta_N^0)^{\cdot}= O\big(\frac{1}{\lambda^3}\big)$
Furthermore, by Proposition \ref{Proposition3.2}
\[R_N(1,\xi_n),S_N(1,\xi_n)= 
\begin{pmatrix}
\ell^2_n\\\ell^2_n
\end{pmatrix}.
\]
Hence taking the $\lambda$--derivative of \eqref{3.30} yields 
\begin{align}
\nonumber
\dot\Delta(\xi_n)= -\dot \theta_N(\xi_n)2\sin \theta_N(\xi_n)+ \frac{1}{(2i\xi_n)^N}\big(\dot R_N(1,\xi_n)_1+ \dot S_N(1,\xi_n)_2\big)+ 
\frac{\ell^2_n}{n^{N+1}}.
\end{align}
By Proposition \ref{Proposition3.4} it then follows that
\begin{align*}
\dot\Delta(\xi_n)=& -\dot \theta_N(\xi_n)2 \sin \theta_N(\xi_n)+ \frac{2i(-1)^{n+1}}{(2i\xi_n)^{N+1}}\int_0^1\varphi_1(t)r_{N+1}(t)\,dt+
\frac{\ell^2_n}{n^{N+1}}\,.
\end{align*}
Going through the arguments of the proof one verifies that the claimed uniformity of the estimate holds.
\end{proof}

\section{Proof of the main results}\label{sec:main_proof}
The aim of this section is to prove the results stated in the introduction. 
\begin{proof}
[Proof of Theorem \ref{Theorem1.5}]
Let $\varphi\in H^N_c$ with $N\geq 1.$ The Dirichlet eigenvalues $\mu_n$ satisfy 
$2i\chi_D(\mu_n) \big(= (m_4+m_3-m_2-m_1)\vert_{(1,\mu_n)}\big)=0.$ By Lemma \ref{Lemma2.2}
in Appendix C, $|\mu_n-n\pi|\leq\frac{1}{|n|}$ for any $|n|\geq n_B.$ Increase $n_B$ if needed so that
$|\alpha_N^0(\mu_n)\beta_N^0(\mu_n)|\leq 1/2$ for any $|n|\geq n_B.$ 
By Proposition \ref{Corollary3.3} (iii) it then follows that 
\begin{align*}
&\frac{\big(1-\alpha_N^{0}(\mu_n)\big)\big(1-\beta_N^{0}(\mu_n)\big)}{1-\alpha_N^{0}(\mu_n)\beta_N^{0}(\mu_n)}  
2 i \sin \theta_N(\mu_n) =\\
& \; i(-1)^{n+1}\big(\hat\varphi_1(-n)+\hat\varphi_2(n)\big)+ \frac{2(-1)^{n}}{(2i\pi n)^{N+1}}\int_0^1\varphi_1r_{N+1}\,dt + 
\frac{\ell^2_n}{n^{N+1}}.
\end{align*}
As $\beta_N^0(\mu_n),\alpha_N^0(\mu_n)=O\big(\frac{1}{n}\big)$ one concludes from the formula above that 
\begin{equation}\label{eq:sine}
\sin \theta_N(\mu_n)= \frac{\ell^2_n}{n^{N}}
\end{equation}
and therefore 
\[
\cos \theta_N(\mu_n)= \pm\sqrt{1-\sin^2\theta_N(\mu_n)}= \pm 1+ \frac{\ell^1_n}{n^{2N}}.
\]
To determine the sign in the above estimate note that $\mu_n=n\pi+O(1/n)$ and by the definition of $\theta_N(\mu_n)$,
\[
\theta_N(\mu_n)= \mu_n+ O\Big(\frac{1}{n}\Big)=n\pi + O\Big(\frac{1}{n}\Big)\,.
\]
Hence, $\cos \theta_N(\mu_n)=(-1)^n+ \frac{\ell^1_n}{n^{2N}}.$ By Proposition \ref{Corollary3.3} (i), it then follows that 
\[
\Delta(\mu_n)= 2(-1)^n+ \frac{\ell^2_n}{n^{N+1}}
\] 
as claimed.  

\noindent As $\frac{\beta_N^{0}(\mu_n)-\alpha_N^{0}(\mu_n)}{1-\alpha_N^{0}(\mu_n)\beta_N^{0}(\mu_n)}=O\big(\frac{1}{n}\big)$ 
one has in view of \eqref{eq:sine} that
\[
\frac{\beta_N^{0}(\mu_n)-\alpha_N^{0}(\mu_n)}{1-\alpha_N^{0}(\mu_n)\beta_N^{0}(\mu_n)}  2i\sin \theta_N(\mu_n)= 
\frac{\ell^2_n}{n^{N+1}}
\]
and thus by Proposition \ref{Corollary3.3} (ii), 
\[
\delta(\mu_n)= i(-1)^n\big(\hat\varphi_1(-n)-\hat\varphi_2(n)\big)+ \frac{\ell^2_n}{n^{N+1}}.
\]
Going through the arguments of the proof one sees that the estimates hold uniformly on bounded sets of potentials $\varphi$ in $H^N_c.$
\end{proof}

\noindent The asymptotics of Theorem \ref{Theorem1.5} can be applied to obtain asymptotics of the eigenvalues of $M(1,\mu_n)$,
referred to as Floquet multipliers of $M(1,\mu_n).$ 
They are given by $\frac{\Delta(\mu_n)\pm \delta(\mu_n)}{2}$ (see e.g. \cite{GK}, p. 50).
By the Wronskian identity their product is $1$ and hence for any $n\in \mathbb{Z}$, $\frac{\Delta(\mu_n)+ \delta(\mu_n)}{2}$
does not vanish. In view of the asymptotics in Theorem \ref{Theorem1.5}, for $|n|$ sufficiently large 
\[
\kappa_n:= 2\log\Big((-1)^n\frac{\Delta(\mu_n)+ \delta(\mu_n)}{2}\Big)
\]
is well defined on bounded sets of $\varphi$'s in $H^1_c$.
\begin{rem1}
Actually, according to \cite{GK}, Theorem 10.3, the $\kappa_n$'s are defined and analytic in a complex neighborhood $W$ of $L^2_r$ in
$L^2_c$ for any $n\in\mathbb{Z}$ and when complemented with the $\mu_n$'s form a system of canonical coordinates on $L^2_r.$ 
\end{rem1}
Theorem \ref{Theorem1.5} leads to the following 
\begin{cor1}\label{Corollary 4.1}
For $\varphi \in H^N_c$ with $N\geq 1,$
\[\kappa_n= i \big(\hat\varphi_1(-n)-\hat\varphi_2(n)\big)+\frac{\ell^2_n}{n^{N+1}} \qquad \text{as}\; |n|\to \infty\]
uniformly on bounded sets of $H^N_c.$
\end{cor1}
\begin{proof}
[Proof of Theorem \ref{Theorem1.1}]
Let $\varphi \in H^N_c$ with $N\geq 1.$ By Lemma \ref{Lemma2.2} in Appendix C, $|\mu_n-n\pi|\leq \frac{1}{|n |}\; \forall |n| \geq n_B.$
Choose $n_B$ bigger if needed so that $|\alpha_N^0(\mu_n)\beta_N^0(\mu_n)|\leq 1/2$ for any $|n|\geq n_B.$  
By Proposition \ref{Corollary3.3} (iii) it then follows that  
\[
2i(-1)^n \sin \theta_N(\mu_n)= 2 T_n
\] 
where 
\[
T_n:= i\,\frac{\hat\varphi_1(-n)+\hat\varphi_2(n)}{2}+ \frac{1}{(2i\mu_n)^{N+1}}\int_0^1\varphi_1r_{N+1}\,dt + 
\frac{\ell^2_n}{n^{N+1}}.
\]
Note that $(-1)^n \sin \theta_N(\mu_n)= \sin\big(\theta_N(\mu_n)-n\pi\big).$ 
Hence $\eta_n:= i (\theta_N(\mu_n)-n\pi)$ satisfies $e^{\eta_n}-e^{-\eta_n}= 2T_n.$ The quadratic equation 
$e^{2\eta_n}-2T_ne^{\eta_n} -1=0$ then yields 
$e^{\eta_n}= T_n +\sqrt[+]{1+T_n^2}.$
By taking the logarithm of both sides of the latter identity and in view of the definition \eqref{3.17bis} and the estimate 
$T_n= \frac{\ell^2_n}{n^N}$ it then follows that 
\begin{align*}
&i (\mu_n-n\pi)- \sum_{k=1}^N\frac{1}{(2i\mu_n)^k}\int_0^1 \varphi_1r_kdt = \eta_n = \\
&  i\,\frac{\hat\varphi_1(-n)+\hat\varphi_2(n)}{2}+ \frac{1}{(2i\mu_n)^{N+1}}\int_0^1\varphi_1r_{N+1}\,dt +
\frac{\ell^2_n}{n^{N+1}}
\end{align*}
leading to 
\begin{align}\label{4.1}
\mu_n=& n\pi-i \sum_{k=1}^{N+1}\frac{1}{(2i\mu_n)^k}\int_0^1 \varphi_1r_k\,dt 
+ \frac{\hat\varphi_1(-n)+\hat\varphi_2(n)}{2}+  \frac{\ell^2_n}{n^{N+1}}.
\end{align}
Unfortunately, $\mu_n$ appears also on the right hand side of the latter asymptotic estimate. 
To address this issue we use an argument applied first by Marchenko in \cite{Ma} (see also, \cite{KST1}, p. 260).
Introduce $F(z):= i \sum_{k=1}^{N+1}\frac{\int_0^1\varphi_1 r_k\,dt}{(2i)^k}\,z^k$ and write
$\zeta_n:=\mu_n-n\pi$ so that
\[
\frac{1}{\mu_n}= \frac{1}{n\pi +\zeta_n}= \frac{\frac{1}{n}}{\pi+\frac{\zeta_n}{n}}.
\]
We approximate $F(\frac{1}{\mu_n})$ by approximating $\zeta_n$ by $\zeta(\frac{1}{n})$ in the above expression for 
$\frac{1}{\mu_n}$ where $\zeta$ is an analytic function so that near $z=0,$
$\zeta(z)+ F(\frac{z}{\pi+ z\zeta(z)})=0.$ 
To find $\zeta$ introduce $G(z,w):= w+F(\frac{z}{\pi+zw})$, defined in an open neighborhood of $(0,0)$ in $\mathbb{C}^2$.
Note that $G$ is analytic, $G(0,0)=0$, and $\partial_w G(0,0)=1.$
Hence by the implicit function theorem there exists near $z=0$ a unique analytic function $\zeta= \zeta(z)$ so that 
$\zeta(0)=0$ and $G(z,\zeta(z))=0$ for $z$ near $0.$
It follows that $\zeta$ has an expansion of the form $\zeta(z)= \sum_{k=1}^\infty c_kz^k.$ 
The coefficients $c_k$ ,$k\geq 1,$ can be computed recursively from the identity  $\zeta(z)=- F(\frac{z}{\pi+ z\zeta(z)}).$ 
In this way one sees that for any $k\geq 1,$ $c_k$ are expressions in $\int_0^1\varphi_1 r_k\,dt$, $1\le k\le N+1$.
Now let us compare $F(\frac{1}{\mu_n})$ with its approximation $F(\frac{1}{\nu_n})$ where $\nu_n:= n\pi + \zeta(\frac{1}{n}).$ 
Using $\frac{1}{\mu_n}-\frac{1}{\nu_n}= \frac{\zeta(\frac{1}{n})-\zeta_n}{\mu_n\nu_n}$ we verify that 
\begin{align}
\label{4.2} 
F\Big(\frac{1}{\mu_n}\Big)= F\Big(\frac{1}{\nu_n}\Big)+ F_n\cdot\Big(\zeta\Big(\frac{1}{n}\Big)-\zeta_n\Big)
\end{align}
where 
\begin{align}\label{4.3}
F_n:= \int_0^1F'\Big(\frac{1}{\nu_n}+t\Big(\frac{1}{\nu_n}-\frac{1}{\mu_n}\Big)\Big)\,dt\cdot \frac{1}{\mu_n\nu_n}
= O\Big(\frac{1}{n^2}\Big)
\end{align} 
as $\zeta(0)=0$ and $\frac{1}{\nu_n}=O\big(\frac{1}{n}\big).$ Rewrite \eqref{4.1} as 
$\zeta_n= -F(\frac{1}{\mu_n})+ \frac{\ell^2_n}{n^N}$ 
and subtract 
\begin{equation}\label{eq:zeta-relation}
\zeta\Big(\frac{1}{n}\Big)= - F\Big(\frac{1}{\nu_n}\Big)
\end{equation}
to get 
\[
\zeta_n-\zeta\Big(\frac{1}{n}\Big)= -\Big(F\Big(\frac{1}{\mu_n}\Big) - F\Big(\frac{1}{\nu_n}\Big)\Big) +\frac{\ell^2_n}{n^N}
\]
implying, in view of \eqref{eq:zeta-relation}, that $(1-F_n)(\zeta_n-\zeta(\frac{1}{n})) = \frac{\ell^2_n}{n^N}.$ 
By \eqref{4.3} one then concludes that $\zeta_n-\zeta(\frac{1}{n})= \frac{\ell^2_n}{n^N}$ which by \eqref{4.2} and \eqref{eq:zeta-relation}
yields 
\[
F\Big(\frac{1}{\mu_n}\Big)+\zeta\Big(\frac{1}{n}\Big)=F\Big(\frac{1}{\mu_n}\Big)-F\Big(\frac{1}{\nu_n}\Big) = \frac{\ell^2_n}{n^{N+2}}.
\]
Altogether we have proved that 
\begin{align*}\mu_n  = & n\pi -F\Big(\frac{1}{\mu_n}\Big)+\frac{\hat\varphi_1(-n)+\hat\varphi_2(n)}{2} +\frac{\ell^2_n}{n^{N+1}}
\\=& n\pi + \zeta\Big(\frac{1}{n}\Big)+\frac{\hat\varphi_1(-n)+\hat\varphi_2(n)}{2} +\frac{1}{n^{N+1}}\ell^2_n
\\=& 
n\pi + \sum_{k=1}^{N+1}\frac{c_k}{n^k}+\frac{\hat\varphi_1(-n)+\hat\varphi_2(n)}{2} +\frac{\ell^2_n}{n^{N+1}}
\end{align*}
which proves the claimed asymptotic estimates. Going through the arguments of the proof one verifies that the estimates hold 
uniformly on bounded sets of $\varphi$'s in $H^N_c.$
\end{proof}

\begin{rem1}\label{rem1}
As mentioned above, the $c_k$'s can be determined recursively from the identity $\zeta(z)= -F(\frac{z}{\pi+z\zeta(z)}).$ One computes 
\[
c_1 = \frac{1}{2\pi}\int_0^1 \varphi_1(t)\varphi_2(t) dt,\qquad 
c_2 = \frac{i}{4\pi^2}\int_0^1 \varphi_1(t)\varphi_2'(t) dt \, .
\]
\end{rem1}

\medskip

\begin{proof}
[Proof of Theorem \ref{Theorem1.2}]
Let $\varphi\in H^N_c$ with $N\geq 1.$ By Lemma \ref{Lemma2.3}, $|\lambda_n^\pm-n\pi|\leq \frac{1}{|n|}$ for $|n|\geq n_B.$ 
Comparing with the case $\varphi= (0,0)$ it then follows that $\Delta(\lambda_{2n}^\pm)=2$ and $\Delta(\lambda_{2n+1}^\pm)=-2$ for $|
2n|\geq n_B.$ It means that $\lambda_{2n}^\pm$ $[\lambda_{2n+1}^\pm]$ are periodic [antiperiodic] eigenvalues of $L(\varphi)$ for 
$|n|\geq \frac{n_B}{2}.$ The proof of the asymptotic estimates of the periodic and antiperiodic eigenvalues are similar so we concentrate
on the asymptotics of the periodic ones only.
Note that the periodic eigenvalues $\lambda_{2n}^\pm,\, |2n|\geq n_B,$ satisfy the equation
\[ 
\det\big(M_N(1,\lambda_{2n}^\pm)-M_N(0,\lambda_{2n}^\pm)\big)=0
\]
By \eqref{3.16}--\eqref{3.17} and Proposition \ref{Proposition3.2}, 
\[M_N(1,\lambda_{2n}^\pm)= \begin{pmatrix}
e^{-i\theta_N}+a_1& \beta_N^0e^{i\theta_N}+ a_2\\
 \alpha_N^0e^{-i\theta_N}+ a_3 & e^{i\theta_N}+a_4
\end{pmatrix}\]
where, with $e_1= \frac{1}{(2i\lambda_{2n}^\pm)^{N+1}}\int_0^1 \varphi_1r_{N+1}dt,$
\[
a_1= e_1+ \frac{1}{n^{N+1}}\ell^2_n,\qquad a_4= -e_1+\frac{1}{n^{N+1}}\ell^2_n
\]
and with $e_2=i \hat \varphi_1(-2n),$ $e_3=- i \hat\varphi_2(2n)$
\[
a_2= e_2  +\frac{1}{n^{N+1}}\ell^2_n,\qquad  a_3= e_3  +\frac{1}{n^{N+1}}\ell^2_n            
\]
and where $\alpha_N^0,\beta_N^0, \theta_N$ are evaluated at $\lambda_{2n}^\pm.$ As 
$M_N(0,\lambda_{2n}^\pm)=\begin{pmatrix}
1& \beta_N^0\\
\alpha_N^0&1
\end{pmatrix}$,
$\det\big(M_N(1,\lambda_{2n}^\pm)-M_N(0,\lambda_{2n}^\pm)\big) $
is given by
\[
(e^{-i\theta_N}-1 +a_1)(e^{i\theta_N}-1 +a_4)- (\alpha_0^Ne^{-i\theta_N}-\alpha_0^N+a_3)(\beta_0^Ne^{i\theta_N}-\beta_0^N+a_2) \, .
\]
Hence $\eta_n\equiv \eta_n^\pm:= e^{i\theta_N(\lambda_{2n}^\pm)}$ satisfies the following quadratic equation\\ 
$a\eta_n^2+b\eta_n +c=0$ where 
\begin{align*}
a= -1+\alpha_N^0\beta_N^0 +a_1-\beta_N^0 a_3, \quad 
c= -1+ \alpha_N^0\beta_N^0 +a_4-\alpha_N^0a_2\\
b= 1+ (1-a_1)(1-a_4) -\alpha_N^0\beta_N^0- (\alpha_N^0-a_3)(\beta_N^0 -a_2)\, .
\end{align*}
Note that $-b=a+c+A$ where $A:= a_2a_3-a_1a_4.$ Hence
$\eta_n= -\frac{b}{2a}+ \frac{1}{2a} \sqrt{b^2-4ac}$ can be written as 
\begin{align}\label{4.5}
\eta_n= \frac{a+c+A}{2a}+ \frac{1}{2a}\sqrt{(a-c)^2 +2A(a+c)+A^2}.
\end{align}
We will address the question of the sign of the root below. First let us analyze the size of the various terms in
the above expression for $\eta_n.$ 
Concerning the term $\frac{a+c+A}{2a}= 1+\frac{c-a+A}{2a},$ note that 
$c-a=-2e_1 + \frac{\ell^2_n}{n^{N+1}},$ $A= \frac{\ell^1_n}{n^{2N}},$ and $2a=-2+O(\frac{1}{n^2})$.
Hence
\begin{align}\label{4.6} 
\frac{c-a+A}{2a}= 1 + e_1 + \frac{\ell^2_n}{n^{N+1}}.
\end{align}
Concerning the expression inside the square root in \eqref{4.5}, one has 
\[(a-c)^2= 4e_1^2+ \frac{\ell^2_n}{n^{2N+2}},\quad  A^2= \frac{\ell^1_n}{n^{4N}},
\] 
and \[
A= e_2e_3+e_1^2+ e_2\frac{\ell^2_n}{n^{N+1}} + e_3\frac{\ell^2_n}{n^{N+1}}
+ \frac{\ell^2_n}{n^{2N+2}}.
\]
As $a+c= -2+O(\frac{1}{n^2})$ one then gets
\begin{align}\label{4.7}
(a-c)^2 +2A(a+c)+A^2= -4e_2e_3 +h^\pm_{2n}
\end{align}
where $e_2e_3= \hat\varphi_1(-2n)\hat\varphi_2(2n)$ and 
\begin{align}\label{4.8}
h_{2n}^\pm= e_2\frac{\ell^2_n}{n^{N+1}}+ e_3 \frac{\ell^2_n}{n^{N+1}} + \frac{\ell^2_n}{n^{2n+2}}= 
\frac{\ell^1_n}{n^{2N+1}}.
\end{align}
Combining these estimates yields
\[
e^{i\theta_N(\lambda_n^\pm)}= 1+e_1 -i\sqrt{e_2e_3+h_{2n}^\pm}+ \frac{\ell^2_n}{n^{N+1}}\, .
\]
Taking the logarithm on both sides of the latter identity then leads to 
\[
\theta_N(\lambda_{2n}^\pm)-2n\pi= -ie_1 + \sqrt{e_2e_3+ h_{2n}^\pm} + \frac{\ell^2_n}{n^{N+1}} \, .
\]
Finally in view of the definition \eqref{3.17bis} of $\theta_N$ we conclude that 
\begin{align}\label{4.9bis}
\lambda_{2n}^\pm= &2n\pi -i \sum_{k=1}^{N+1}\frac{1}{(2i\lambda_{2n}^\pm)^k}\int_0^1\varphi_1r_k\,dt 
+ \sqrt{e_2e_3+ h_{2n}^\pm} + \frac{\ell^2_n}{n^{N+1}}.
\end{align}
To address the issue of the sign of the root in \eqref{4.5}, introduce
\[
A^\pm:= \big\{|2n|\geq n_B \,\,\big|\,\, |h_{2n}^\pm|<|e_2e_3|/2\big\}.
\]
It then follows that for $|2n| \ge n_B$ with $2n \notin A^+\cap A^-,$ 
$|e_2e_3 + h_{2n}^\pm|^{\frac{1}{2}}= \frac{\ell^2_n}{|n|^{N+\frac{1}{2}}}$
implying that $|e_2e_3|^{\frac{1}{2}}= \frac{\ell^2_n}{|n|^{N+\frac{1}{2}}}.$
For $2n\in A^+\cap A^-,$ note that $e_2e_3\neq 0.$ Denote by $\sqrt[\circ ]{e_2e_3}$ an arbitrary branch of the square root 
(which might depend on $n$) and by $\sigma_n^\pm\in \{1,-1\}$ the sign of the root determined by \eqref{4.5} so that
\begin{align*}
\lambda_{2n}^\pm= 2n\pi-i \sum_{k=1}^{N+1}\frac{1}{(2i\lambda_{2n}^\pm)^k}\int_0^1\varphi_1r_kdt + 
\sigma_n^\pm\sqrt[\circ ]{e_2e_3+h_{2n}^\pm}+\frac{\ell^2_n}{n^{N+1}}.
\end{align*} 
Let $A_0:= \{2n\in A^+\cap A^-|\; \sigma_n^+=\sigma_n^-\}.$
(Note that $A_0$ could be empty or finite). By Lemma \ref{LemmaB.2}, 
$\Big|\sqrt[\circ ]{e_2e_3 +h_{2n}^\pm}-\sqrt[\circ ]{e_2e_3}\Big|\leq |h_{2n}^\pm|^{1/2}$ for any $2n\in A_0.$
As by Theorem \ref{Theorem1.4}, for any $2n\in A_0$
\[
\sigma_n^+\sqrt[\circ ]{e_2e_3+h_{2n}^+}+\sigma_n^-\sqrt[\circ ]{e_2e_3+h_{2n}^-}= \frac{\ell^2_n }{n^{N+1}}
\]
it follows from Lemma \ref{LemmaB.2} that 
\begin{align*}
|2\sqrt[\circ ]{e_2e_3}|\leq \, \, &\Big|\sqrt[\circ ]{e_2e_3+h_{2n}^+}+ \sqrt[\circ ]{e_2e_3+h_{2n}^-}\Big| \\\, \, &+
\Big|\sqrt[\circ ]{e_2e_3+h_{2n}^+}- \sqrt[\circ ]{e_2e_3}\Big|+
\Big|\sqrt[\circ ]{e_2e_3+h_{2n}^-}- \sqrt[\circ ]{e_2e_3}\Big| 
\\ =\, \,&
\frac{\ell^2_n}{n^{N+1}} + |h_{2n}^+|^{1/2}+|h_{2n}^-|^{1/2}= \frac{\ell^2_n}{n^{N+\frac{1}{2}}}.
\end{align*}
Hence we have proved the following asymptotic estimates 
\begin{equation*}
\{\lambda^+_{2n},\lambda^-_{2n}\}= \Big\{ 
2n\pi -i\sum_{k=1}^{N+1}\frac{1}{(2i\lambda_{2n}^\pm)^k}\int_0^1 \varphi_1r_k\,dt \pm\sqrt{\hat \varphi_1(-2n) \hat\varphi_2(2n)}
+\frac{\ell^2_n}{n^{N+\frac{1}{2}}}\Big\}\,.
\end{equation*}
By applying as in the proof of Theorem \ref{Theorem1.1} Marchenko's argument one obtains the claimed asymptotics of item (i). 

\noindent Towards item (ii) we first remark that for $\varphi\in H^N_r,$ $\hat\varphi_1(-2n)= \overline{\hat\varphi_2}(2n)$ and therefore
$e_2e_3= |\hat\varphi_1(-2n)|^2$.
Our starting point is formula \eqref{4.9bis}. As in the case at hand $|e_2|= |e_3|$ we can write $e_2e_3+ h^\pm_{2n},$ 
given by \eqref{4.8} as follows
\[
e_2 e_3+h_{2n}^\pm= (|e_2|+g_{2n}^\pm)^2+k^\pm_{2n}
\]
where $g_{2n}^\pm= \frac{\ell^2_n}{n^{N+1}}$ and  $k^\pm_{2n}=\frac{\ell^2_n}{n^{2N+2}}.$ 
Now define 
\[
A^\pm= \Big\{|2n|\geq n_B \,\,\Big|\,\, |k^\pm_{2n}|\leq \frac{ | (|e_2| + g^\pm_{2n})|^2 }{2}\Big\}.
\]
For $|2n| \ge n_B$ with $2n \notin  A^+,$
$|(|e_2|+g^+_{2n})|^2= \frac{\ell^2_n}{n^{2N+2}},$ implying that 
$|e_2|+g^+_{2n}=\frac{\ell^4_n}{n^{N+1}}$ and hence $|e_2|= \frac{\ell^4_n}{n^{N+1}}.$
Similarly, if $|2n| \ge n_B$ with $2n \notin  A^-,$ $|e_2|= \frac{\ell^4_n}{n^{N+1}}.$ 
If $2n\in A^+\cap A^-,$ then by Lemma \ref{LemmaB.2} (i), 
\[\sqrt[\circ ]{(|e_2|+g^\pm_{2n})^2 + k_{2n}^\pm}=\sqrt[\circ ]{(|e_2|+g^\pm_{2n})^2 }+ \frac{\ell^4_n}{n^{N+1}}\]
where $\sqrt[\circ ]{\cdot}$ denotes an arbitrary branch of the square root. 
Arguing as in the proof of item (i) and taking into account that $\lambda_{2n}^-\leq \lambda_{2n}^+$ the claimed asymptotics 
\[
\lambda_{2n}^\pm =
2n\pi -i\sum_{k=1}^{N+1}\frac{1}{(2i\lambda_{2n}^\pm)^k}\int_0^1 \varphi_1r_k dt \pm|\hat \varphi_1(-2n)| +
\frac{\ell^4_n }{n^{N+1}}
\]
follow. Going through the arguments of the proofs of (i) and (ii) one verifies that the stated uniformity property holds.
\end{proof}

\begin{proof}
[Proof of Theorem \ref{Theorem1.4} (i)]
By Lemma \ref{Lemma2.2}, $|\dot \lambda_n- n\pi|\leq \frac{1}{|n|},$ for any $|n|\geq n_B.$ 
By Corollary \ref{Corollary3.5} for  $|n|\geq n_B$ 
\begin{align}\label{4.20}
\dot \Delta(\dot \lambda_n) = -\dot\theta_N2 \sin \theta_N\big|_{\lambda = \dot \lambda_n}+ 2i 
\frac{(-1)^{n+1}}{(2i \dot \lambda_n)^{N+1}}\int_0^1\varphi_1 r_{N+1}\,dt 
+ \frac{\ell^2_n}{n^{N+1}} \, .
\end{align}
By \eqref{3.17bis}, $\theta_N(\lambda)= \lambda + i\sum_{k=1}^N\frac{1}{(2i\lambda)^k}\int_0^1 \varphi_1r_k\,dt$ and 
hence $\dot\theta_N(\dot \lambda_n)= 1+O(\frac{1}{n^2}).$ Therefore $\dot\Delta_N(\dot \lambda_n)=0$ yields
\[
\sin\theta_N(\dot \lambda_n)+ i a_n=0
\] 
where $a_n= \frac{(-1)^n}{(2i\dot\lambda_n)^{N+1}}\int_0^1\varphi_1r_{N+1}\,dt+ \frac{\ell^2_n}{n^{N+1}}.$ 
Introduce the following sequence $\eta_n:= e^{i\theta_N(\dot\lambda_n)}.$
As $2i \sin\theta_N(\dot\lambda_n)= \eta_n-\eta_n^{-1}$ it then follows that $\eta_n^2- 2a_n\eta_n-1=0$
implying that $\eta_n= a_n+(-1)^n\sqrt[+]{1+a_n^2}= (-1)^n+a_n +O(a_n^2).$
Taking the logarithm on both sides of the latter identity leads to 
$\theta_N(\dot \lambda_n)= n\pi -i(-1)^na_n+O(a_n^2)$ or 
\[
\dot\lambda_n= n\pi -i\sum_{k=1}^{N+1}\frac{1}{(2i\dot\lambda_n)^k}\int_0^1 \varphi_1r_k\,dt+ \frac{\ell^2_n}{n^{N+1}}.
\]
Arguing as in the proof of Theorem \ref{Theorem1.1} (use Marchenko's argument) it follows that 
 \[
\dot\lambda_n= n\pi -i\sum_{k=1}^{N+1}\frac{c_k}{n^k}+ \frac{\ell^2_n}{n^{N+1}}.
\]
Going through the arguments of the proof one verifies that the latter estimate holds uniformly on bounded sets of $H^N_c.$
\end{proof}

\begin{proof}
[Proof of Theorem \ref{Theorem1.1} (ii)]
According to \cite{GK}, Lemma 6.9, for any $\varphi\in L^2_c$ there exists $n_0\geq 1$ and a neighborhood 
$V$ of $\varphi$ in $L^2_c$ so that 
\[
\tau_n= \dot \lambda_n+ O(\gamma_n^2)\quad \forall |n|\geq n_0
\] 
uniformly on $V.$ As $H^1_c \hookrightarrow L^2_c$ is a compact embedding it follows that 
$\tau_n= \dot \lambda_n+ O(\gamma_n^2)$ uniformly on bounded sets of $H^N_c$ with $N\geq 1.$
The claimed asymptotics of $\tau_n$ then follow from item (i) of Theorem \ref{Theorem1.4} and Corollary \ref{Corollary1.3}. 
\end{proof}

\section{Appendix A: Asymptotic estimates of M}
In this appendix we prove asymptotic estimates of the fundamental solution $M(x,\lambda)$ of the linear system 
$L(\varphi) M=\lambda M$ for $\varphi\in H^1_c$. 
Recall (see e.g. \cite{GK}, Section 1) that for $\varphi \in L^2_c$, $M \equiv M(x, \lambda)$ is a continuous function on 
$[0,1]\times\mathbb C$, given by the infinite series $M= \sum_{n=0}^\infty M_n$ with 
$M_0(x,\lambda)=E_\lambda(x) $ and, for any $n\geq 0,$ 
\[
M_{n+1}(x,\lambda) = \int_0^xE_{\lambda}(x-x_1)R \Phi(x_1) M_n(x_1,\lambda)dx_1
\] 
where
\[
E_\lambda(x)= \begin{pmatrix}
e^{-i\lambda x}& 0\\ 0 &e^{i\lambda x}
\end{pmatrix},\quad R= \begin{pmatrix}
i &0 \\ 0& -i
\end{pmatrix},\quad \Phi (x) = \begin{pmatrix}
0 & \varphi_1(x) \\ \varphi_2(x) & 0 
\end{pmatrix}  \, .
\]
Note that for any $n\geq 0,$ $M_{2n}$ is a diagonal $2\times 2$ matrix whereas $M_{2n+1}$ is off-diagonal. 
In the sequel we will always assume that $\varphi\in H^1_c$ if not stated otherwise. 
Then $M(x, \lambda)$ is a continuously differentiable function in $0 \le x \le 1$ and 
$\lambda \in \mathbb C$.
Throughout the appendix we will use the elementary identities
\begin{align}\label{A.0}
\Phi(x)E_\lambda(x)= E_\lambda(-x)\Phi(x),\, [R,E_\lambda(x)]=0, \,R^2=-1, \,\,\,\text{and}\,\,R\Phi=-\Phi R.
\end{align}
We begin by taking a closer look at $M_1(x,\lambda), M_2(x,\lambda)$, and $M_3(x,\lambda).$ 
By \eqref{A.0}, one has $M_1(x,\lambda)=\int_0^xE_\lambda(x-2t)R\Phi(t)\,dt.$ 
Integrating by parts and taking into account that for $\lambda\in\mathbb{C}\setminus\{0\}$,
\begin{align}\label{A.0bis}
E_\lambda(x-2t)= -\frac{1}{2\lambda}\,R\,\partial_t \big(E_\lambda(x-2t)\big)
\end{align}
we get
\begin{align}\label{A.1}
M_1(x,\lambda)= \frac{1}{2\lambda }\,\big(E_\lambda(-x)\Phi(x)- E_\lambda(x)\Phi(0) - E_\lambda(x)P_\lambda(x)\big)
\end{align}
where 
\begin{align}\label{A.1bis}
P_\lambda(x):= \int_0^x E_\lambda(-2t)\Phi'(t)\,dt, \quad \Phi'(t)=\partial_t\Phi(t).
\end{align}
Substituting the expression \eqref{A.1} for $M_1$ into the expression 
\[
M_{2}(x,\lambda) = \int_0^xE_{\lambda}(x-x_1) R \Phi(x_1) M_1(x_1,\lambda)\,dx_1
\] 
one gets $M_2(x,\lambda) = \frac{1}{2\lambda}(I+II+III)$ where
\[
I:=\int_0^xE_{\lambda}(x-x_1) R \Phi(x_1) E_\lambda(-x_1)\Phi(x_1)\,dx_1 
\]
leading in view of \eqref{A.0} to 
\begin{align}
\label{A.2}
I=& E_\lambda(x) R Q(x),\qquad Q(x):= \int_0^x \varphi_1(t)\varphi_2(t)\,dt,
\\\label{A.3}
II:= & -\int_0^xE_{\lambda}(x-x_1)R \Phi(x_1) E_\lambda(x_1)\Phi(0)\,dx_1 = -M_1(x,\lambda)\Phi(0)
\end{align}
and $III:= -\int_0^xE_{\lambda}(x-2x_1)R \Phi(x_1) P_\lambda(x_1)\,dx_1.$
The latter term can be integrated by parts to get with \eqref{A.0bis}
\begin{align}
III= -\frac{1}{2\lambda}\int_0^x \big(E_\lambda(x-2x_1)\Phi(x_1)P_\lambda(x_1)\big)'\,dx_1+&\\
+\frac{1}{2\lambda}\int_0^x E_\lambda(x-2x_1)\big(\Phi(x_1)P_\lambda(x_1)\big)'\,dx_1.
\end{align}
As $P_\lambda(0)=0$ and $\Phi(x_1)E_{\lambda}(-2x_1)= E_{\lambda}(2x_1)\Phi(x_1)$ one gets
\begin{align}\nonumber 
III=& -\frac{1}{2\lambda}E_{\lambda}(-x)\Phi(x)P_\lambda(x)+ \frac{1}{2\lambda}E_{\lambda}(x)\int_0^x \Phi(x_1)\Phi'(x_1)\,dx_1 +
\\\label{A.4}&+\frac{1}{2\lambda}\int_0^xE_{\lambda}(x-2x_1)\Phi'(x_1)P_\lambda(x_1)\,dx_1 \, .
\end{align} 
Combining \eqref{A.2}-\eqref{A.4} then yields
\begin{align}\label{A.5}
M_2(x,\lambda)-\frac{1}{2\lambda}E_\lambda(x) RQ(x)= - \frac{1}{4\lambda^2}\hat{M_2}(x,\lambda)
\end{align}
where 
\begin{align}
\nonumber \hat{M_2}(x,\lambda)=& 2\lambda M_1(x,\lambda)\Phi(0)+ E_\lambda(-x)\Phi(x)P_\lambda(x) 
-E_\lambda(x)\int_0^x\Phi(x_1)\Phi'(x_1)dx_1 \\\label{A.6} & - 
\int_0^xE_{\lambda}(x-2x_1)\Phi'(x_1)P_\lambda(x_1)\,dx_1.
\end{align}
where, in view of \eqref{A.1},
\begin{equation}\label{eq:app_rel1}
2\lambda M_1(x,\lambda)\Phi(0)=E_\lambda(-x)\Phi(x)\Phi(0)-E_\lambda(x)\Phi(0)^2-
E_\lambda(x)P_{\lambda}(x)\Phi(0)\,.
\end{equation}
Let $|A|$ be the operator norm $|A|:=\max\limits_{|x|=1}|Ax|$ of a complex matrix, $A:=(a_{kl})_{1\le k,l\le 2}$,
where $|x|=\sqrt{|x_1|^2+|x_2|^2}$, $x\in\mathbb{C}^2$. 
Note that for any $a, b\in\mathbb{C},$ 
\[
\left|\begin{pmatrix}
a&0 \\ 0&b
\end{pmatrix}\right|=\left|\begin{pmatrix}
0&a \\ b&0
\end{pmatrix}\right|= \max{(|a|,|b|)}\,.
\]
One easily sees that for any $y\in\mathbb{R}$,
\begin{equation}\label{eq:app_rel2}
|E_\lambda(y)|\le e^{|\operatorname{Im}\lambda| |y|}\,.
\end{equation}
In particular, for any $0\leq t\leq x \leq 1,$ one has 
$|E_\lambda(x-2t)|\leq e^{ |\operatorname{Im}\lambda |x}$,
\[
 |E_\lambda(x)P_\lambda(x)|\leq  \int_0^x|E_\lambda(x-2t)||\Phi'(t)|\,dt
\leq e^{|\operatorname{Im}\lambda|x}\|\varphi\|_{H^1},
\]
where we used that $\max(|a|, |b|)\le |a|+|b|$.
Note that by Sobolev embedding, 
\begin{align}\label{A.7}
|\Phi(x)|= \max_{i=1,2}|\varphi_i(x)|\leq c\,\|\varphi\|_{H^1}
\end{align}
for some constant $c>0$.
Using formula \eqref{eq:app_rel1} for $2\lambda M_1(x,\lambda)$ one verifies that for an absolute constant $C>0$
\[
|2\lambda M_1(x,\lambda)\Phi(0)|\leq C e^{ |\operatorname{Im}\lambda |x}\|\varphi\|_{H^1}^2.
\]
The other terms in the formula \eqref{A.6} for $\hat M_2(x,\lambda)$ are estimated in a similar way, yielding 
$|E_\lambda(-x)\Phi(x)P_\lambda(x)|\leq c\,e^{|\operatorname{Im}\lambda|x}\|\varphi\|_{H^1}^2,$
\[ 
\Big|E_\lambda(x)\int_0^x\Phi(x_1)\Phi'(x_1)\,dx_1\Big|
\leq e^{|\operatorname{Im}\lambda|x}\|\varphi\|_{H^1}^2,
\]
and, taking into account \eqref{eq:app_rel2} and
\[
|x-2x_1+2x_2|= |(x-x_1)-(x_1-x_2)+x_2|\leq |x-x_1|+|x_1-x_2|+ x_2= x,
\]
for any $0\leq x_2\leq x_1\leq x,$ 
\begin{align*}
\Big|\int_0^x &E_{\lambda}(x-2x_1)\Phi'(x_1)P_\lambda(x_1)\,dx_1\Big| \\ & \leq
\int_0^x|\Phi'(x_1)|\Big(\int_0^{x_1}|E_{\lambda}(-x+2x_1-2x_2)||\Phi'(x_2)|dx_2\Big)\,dx_1 \\
& \leq  e^{|\operatorname{Im}\lambda|x} \int_0^x|\Phi'(x_1)|\Big(\int_0^{x_1}|\Phi'(x_2)|\,dx_2\Big)\,dx_1 
\leq e^{|\operatorname{Im}\lambda|x}\|\varphi\|_{H^1}^2.
\end{align*}
We thus have proved that there exists an absolute constant $C>0$ so that 
\begin{align}\label{A.8}
|\hat M_2(x,\lambda)|\leq Ce^{|\operatorname{Im}\lambda|x}\|\varphi\|_{H^1}^2,\qquad 
\forall 0\leq x\leq 1,\,\lambda\in \mathbb{C},\,\varphi\in H^1_c.
\end{align}
Next we consider $M_3(x,\lambda)= \int_0^x E_\lambda(x-x_1)R\Phi(x_1)M_2(x_1,\lambda)dx_1.$ Note that by \eqref{A.5},
\[M_3(x,\lambda)=\frac{1}{2\lambda }IV-\frac{1}{4\lambda^2}\int_0^x E_\lambda(x-x_1)R \Phi(x_1) \hat M_2(x_1,\lambda)\,dx_1
\]
where 
\[
IV := \int_0^x E_\lambda(x-x_1)R \Phi(x_1) E_\lambda(x_1)RQ(x_1)\,dx_1 .
\]
One concludes form \eqref{A.8} that for any $\lambda \in \mathbb{C}\setminus \{0\},$ $0\leq x\leq 1,\;\varphi \in H^1_c$
\begin{align*}
|M_3(x,\lambda)-\frac{1}{2\lambda}IV|\leq & \frac{1}{4|\lambda|^2}\int_0^x |E_\lambda(x-x_1)R \Phi(x_1) \hat M_2(x_1,\lambda)|\,dx_1 \\
\leq & \frac{C}{4|\lambda|^2}e^{|\operatorname{Im}\lambda|x}\|\varphi\|_{H^1}^3
\end{align*}
where we used that $|E_\lambda(x-x_1)|\le e^{|\operatorname{Im}\lambda|(x-x_1)}$, $0\le x_1\le x$.
The term $IV$ can be integrated by parts. As $R^2=-Id_{2\times 2}$ and $\Phi R= -R\Phi$ one gets
$IV= \int_0^xE_\lambda(x-2x_1)\Phi(x_1)Q(x_1)\,dx_1.$ As by \eqref{A.0bis},
integration by parts then yields
\[
IV= -\frac{1}{2\lambda}\int_0^x\partial_{x_1}\big(E_\lambda(x-2x_1)\big) R \Phi(x_1)Q(x_1)\,dx_1= 
-\frac{1}{2\lambda}V+\frac{1}{2\lambda}VI
\]
where
\[
V :=E_\lambda(x-2x_1)\Phi(x_1)Q(x_1)\big\vert_{x_1=0}^x,\,
VI :=\int_0^xE_\lambda(x-2x_1) R \big(\Phi(x_1)Q(x_1)\big)'\,dx_1 .
\]
As $Q(0)=0,$ we get in view of \eqref{A.7} that
\[
|V|= \Big|E_\lambda(-x)\Phi(x)\int_0^x \varphi_1(t)\varphi_2(t)\,dt\Big|\leq 
c e^{|\operatorname{Im}\lambda|x}\|\varphi\|_{H^1}^3
\]
and
\begin{align*}
| VI |\leq & \int_0^x|E_\lambda(x-2x_1) R\Phi'(x_1)|\Big(\int_0^x |\varphi_1(t)\varphi_2(t)|\,dt \Big)\,dx_1 \\
&+ \int_0^x|E_\lambda(x-2x_1) R\Phi(x_1)| |\varphi_1(x_1)\varphi_2(x_1)|dx_1 
\\\leq &\; e^{|\operatorname{Im}\lambda|x}(1+c)\|\varphi\|_{H^1}^3.
\end{align*}
Altogether one then gets 
$| \frac{1}{2\lambda}IV | \leq \frac{1}{4|\lambda|^2}e^{|\operatorname{Im}\lambda|x}(1+2c)\|\varphi\|_{H^1}^3$
hence for any $\lambda\in \mathbb{C}\setminus \{0\},\,0\leq x\leq 1,$ and $\varphi\in H^1_c$
\begin{align}\label{A.9}
|M_3(x,\lambda)|\leq \frac{C_1}{4|\lambda|^2}e^{|\operatorname{Im}\lambda|x}\|\varphi\|_{H^1}^3
\end{align}
where $C_1= 1+2c+C$. 
Finally, for any $n\geq 1,$ $M_{n+3}(x,\lambda)$ can be written as
\[
\int\limits_{0\leq x_n\leq \dots \leq x_1\leq x }\!\!\!\!\! E_\lambda\Big(x+2\sum_{k=1}^{n-1}(-1)^kx_k+ (-1)^nx_n\Big) 
\big( \prod_{j=1}^nR \Phi(x_j) \big) M_3(x_n,\lambda) dx_n\cdots dx_1.
\]
Similarly as above, one has for any sequence $0\leq x_n\leq \cdots \leq x_1\leq x$,
\[\Big|x+2\sum_{k=1}^{n-1}(-1)^k x_k+ (-1)^n x_n\Big|\leq |x-x_1|+ |x_1-x_2|+\cdots + |x_{n-1}-x_n|= x-x_n.\]
Hence 
$|E_\lambda\big(x+2\sum_{k=1}^{n-1}(-1)^kx_k+ (-1)^nx_n\big)|
\leq e^{|\operatorname{Im}\lambda|(x-x_n)}.$ With \eqref{A.9} it then follows that
\begin{align*}
|M_{n+3}(x,\lambda)|\leq & \frac{C_1}{4|\lambda|^2}e^{|\operatorname{Im}\lambda|x}\|\varphi\|_{H^1}^3
\frac{1}{n!}\Big(\int_0^x|\Phi(t)|\,dt\Big)^n
\\
\leq &\frac{C_1}{4|\lambda|^2}e^{|\operatorname{Im}\lambda|x}\|\varphi\|_{H^1}^3\frac{1}{n!} \|\varphi\|_{L^2}^n\,.
\end{align*}
Combining this with \eqref{A.5} and \eqref{A.8} we get the following estimate for $M= \sum_{n=0}^\infty M_n:$
\begin{thm1}\label{TheoremA.1}
There exists an absolute constant $C>0$ so that for any $0\leq x\leq 1,$ $\lambda\in\mathbb{C}\setminus \{0\},$ and $\varphi\in H^1_c$,
\begin{align*} &
\Big| M(x,\lambda) - E_\lambda(x)- M_1(x,\lambda) -\frac{1}{2\lambda}E_\lambda(x)R\int_0^x \varphi_1(t)\varphi_2(t)\,dt\Big|\\
& \leq \frac{C}{|\lambda|^2} 
e^{|\operatorname{Im}\lambda|x}e^{\|\varphi\|_{L^2}}\|\varphi\|_{H^1}^2(1+\|\varphi\|_{H^1})
\end{align*}
where $M_1(x,\lambda)= \int_0^x E_\lambda(x-2t)R\Phi(t)\,dt$ equals
\[
\frac{1}{2\lambda}\bigg(E_\lambda(-x)\Phi(x)- E_\lambda(x)\Phi(0) - \int_0^x E_\lambda(x-2t) \Phi'(t) dt\bigg).\]
\end{thm1}
\begin{rem1}
Note that  $E_\lambda(x)+ M_1(x,\lambda) +\frac{1}{2\lambda}E_\lambda(x)R\int_0^x \varphi_1(t)\varphi_2(t)\,dt$
is an approximation of $M(x,\lambda)$ for $|\lambda|$ large where
$ M_1(x,\lambda)$ is off-diagonal and
 $ E_\lambda(x)+\frac{1}{2\lambda}E_\lambda(x)R\int_0^x \varphi_1(t)\varphi_2(t)\,dt$ is a diagonal matrix.
\end{rem1}
\noindent Theorem \ref{TheoremA.1} leads to similar estimates for the inverse of $M(x,\lambda).$ 
By the Wronskian identity, $\det M(x,\lambda) =1,$ the inverse of 
$
M= \begin{pmatrix}
m_1 &m_2 \\m_3 &m_4
\end{pmatrix}
$ 
is given by 
\begin{equation}\label{eq:app_rel3}
M^{-1}= M^{\vee}:=\begin{pmatrix}
m_4& -m_2 \\
-m_3 & m_1
\end{pmatrix}.
\end{equation}
As for any $2\times 2$ matrix $A=\begin{pmatrix}
a_1& a_2\\ a_3 & a_4
\end{pmatrix} $ one has 
\[
|a_j|\leq |A| \leq |A|_\infty:=|a_1|+|a_2|+|a_3|+|a_4|\quad\forall 1\leq j \leq 4
\]
we get that
\begin{equation}\label{eq:app_rel4}
|M^\vee|\le |M^\vee|_\infty=|M|_\infty\le 4|M|.
\end{equation}
The asymptotics of Theorem \ref{TheoremA.1} together with \eqref{eq:app_rel3} and \eqref{eq:app_rel4} lead 
to the following asymptotics of $M(x,\lambda)^{-1}$.
\begin{cor1}\label{CorollaryA.2}
For any $0\leq x\leq 1,$ $\lambda\in\mathbb{C}\setminus \{0\},$ and $\varphi\in H^1_c$ 
\begin{align*} &
\Big| M(x,\lambda)^{-1} - E_\lambda(-x)+ M_1(x,\lambda) +\frac{1}{2\lambda}E_\lambda(-x)R\int_0^x \varphi_1(t)\varphi_2(t)\,dt\Big|\\
& \leq  \frac{4C}{|\lambda|^2} e^{|\operatorname{Im}\lambda|x}
e^{\|\varphi\|_{L^2}}\|\varphi\|_{H^1}^2(1+\|\varphi\|_{H^1})
\end{align*}
where $C>0$ is the same constant as in Theorem \ref{TheoremA.1}.
\end{cor1}
\noindent Theorem \ref{TheoremA.1} and Corollary \ref{CorollaryA.2} are used to obtain asymptotic estimates for the solution of 
the inhomogeneous equation $(L(\varphi) -\lambda)F= f$
\begin{align}\label{A.10}
\big(L(\varphi) -\lambda\big)F= f,\quad F(0,\lambda) = \big(F_1(0,\lambda),F_2(0,\lambda)\big)=(0,0)
\end{align} 
where $f= (f_1,f_2)\in L^2_c,$ $\lambda\in\mathbb{C}\setminus \{0\},$ and $\varphi\in H^1_c.$
Substitute the ansatz $F(x,\lambda) = M(x,\lambda) c(x,\lambda)$ of the method of the variation of parameters into 
equation \eqref{A.10} and use that $R^{-1}=-R$ to see that 
\begin{align}\label{A.12}
F(x,\lambda)= -M(x,\lambda)\int_0^x M(t,\lambda)^{-1} R f(t)\,dt.
\end{align}
By Theorem \ref{TheoremA.1} and Corollary \ref{CorollaryA.2}, 
\begin{align*}
F(x,\lambda)=& -\Big(E_\lambda(x)+ M_1(x,\lambda) +\frac{1}{2\lambda}E_\lambda(x) R Q(x) +O\Big(\frac{1}{\lambda^2}\Big)\Big)
\\ &\int_0^x \Big(E_\lambda(-t)-M_1(t,\lambda) -\frac{1}{2\lambda}E_\lambda(-t) R Q(t) +
O\Big(\frac{1}{\lambda^2}\Big)\Big) R f(t)\,dt
\end{align*}  
leading to the following 
\begin{cor1}\label{CorollaryA.3}
For any $\varphi\in H^1_c,$ and $f\in L^2_c,$
the solution $F(x,\lambda)$ of \eqref{A.10} admits for $|\lambda| \to \infty$ the asymptotic expansion 
\[
F(x,\lambda)= A(x,\lambda) + \frac{1}{2\lambda}\sum_{k=1}^4B_k(x,\lambda)+ O\Big(\frac{1}{\lambda^2}\|f\|_{L^2}\Big)
\]
where $A(x,\lambda)= -\int_0^xR E_\lambda(x-t)f(t)\,dt ,$
\[
B_1(x,\lambda)= -2\lambda M_1(x,\lambda)R\int_0^x E_\lambda(-t) f(t)\,dt,
\]
\[
B_2(x,\lambda)= \; Q(x) \int_0^x E_\lambda(x-t)f(t)\,dt,
\]
\[
B_3(x,\lambda)= \;E_\lambda(x) \int_0^x 2\lambda M_1(t,\lambda)Rf(t)\,dt,\,
\]
\[
B_4(x,\lambda)= -\int_0^x E_\lambda(x-t)Q(t)f(t)\,dt,
\]
with $M_1(x,\lambda)$ as defined in \eqref{A.1} and $Q(x)$ as in \eqref{A.2}\,.
For any $B>0$ and $\Lambda_{\rm Im}>0$ the estimate above is uniform in $0\leq x \leq 1$, 
$|\operatorname{\rm Im}\lambda|\le\Lambda_{\rm Im}$, and $\|\varphi\|_{H^1}\le B$.
\end{cor1}

\section{Appendix B: Auxilary lemmas}
First we prove an estimate on perturbed Fourier coefficients used throughout the paper. 
See e.g. \cite{GK}, Appendix D, for similar results and references.
\begin{lem1}\label{LemmaB.1} Let $f\in L^2([0,1], \mathbb{C})$ and let 
\[\phi_n(x)= \int_0^x e^{i\xi_n(x-2t)}f(t) dt, \quad n\in \mathbb{Z}\]
with a sequence of complex numbers $\xi_n= n\pi +\alpha_n$ such that $|\alpha_n|\leq \frac{a}{\langle n\rangle} $ for any 
$n\in \mathbb{Z}$ with $a>0,\,$ and  $\langle n\rangle= \max(1,|n|)$. 
Then for any $0\leq x\leq 1,$
\[\sum_{n\in \mathbb{Z}}\langle n\rangle^2\left|\phi_n(x)-\int_0^x e^{i\pi n(x-2t)}f(t)\,dt\right|^2\leq e^{2a}\|f\|^2_{L^2}.\]
In particular, for $x=1,\; \int_0^1 e^{i\pi n(1-2t)}f(t)\,dt = (-1)^n\hat f(n)$ and hence 
$$
\Big(\sum_{n\in \mathbb{Z}}\langle n\rangle^2\big|\phi_n(1)-(-1)^n\hat f(n)\Big|^2\big)^{1/2}\leq e^{a}\|f\|_{L^2}.
$$
\end{lem1}
\begin{proof}
Setting $g_n(x)=\phi_n(x)-\int_0^x e^{i\pi n(x-2t)}f(t)\,dt $ one gets 
\[
g_n(x)= \int_0^x e^{i\pi n(x-2t)}\left(e^{i\alpha_n(x-2t)}-1\right)f(t)\,dt \, .
\]
Expanding $e^{i\alpha_n(x-2t)}$ into a power series in $x-2t,$ one obtains 
\[g_n(x)= \sum_{k=1}^\infty \frac{(i\alpha_n)^k}{k!}\int_0^{x}(x-2t)^kf(t)e^{i\pi n(x-2t)}\,dt.\]
Denoting the last integral by $f_{k,x,n}$ and using that $|\alpha_n|\leq a\frac{1}{\langle n\rangle} $ yields
\begin{align}
\label{B.1}
|g_n(x)|\leq\frac{1}{\langle n\rangle} \sum_{k=1}^\infty \frac{a^k}{k!}|f_{k,x,n}| \, .
\end{align}
Note that $f_{k,x,n}$ is the n'th Fourier coefficient $\hat f_{k,x}(n)$ of  
$$
f_{k,x}(t):= (x-2t)^kf(t)e^{i\pi n x}\mathbbm{1}_{[0,x]}(t)\, .
$$
Multiplying \eqref{B.1} with $\langle n\rangle^2|g_n(x)|,$ we get
\[
\sum_{n\in \mathbb{Z}}\langle n\rangle^2|g_n(x)|^2 \leq 
\sum_{k=1}^\infty \frac{a^k}{k!}\sum_{n\in \mathbb{Z}}|\hat f_{k,x}(n)|\langle n\rangle|g_n (x)|
\]
which by Cauchy--Schwarz can be bounded by
\[
\sum_{k=1}^\infty \frac{a^k}{k!}\big(\sum_{n\in \mathbb{Z}}\langle n\rangle^2|g_n (x)|^2\big)^\frac{1}{2}
\big( \sum_{n\in \mathbb{Z}}|\hat f_{k,x}(n)|^2\big)^\frac{1}{2}
\leq e^a\big(\sum_{n\in \mathbb{Z}}\langle n\rangle^2|g_n (x)|^2\big)^{\frac{1}{2}}\|f_{k,x}\|_{L^2}\, ,
\]
implying the claimed estimate.
\end{proof}
\noindent In Section 4 we frequently use the following elementary estimates of the square root. For $z\in \mathbb{C}\setminus \{0\}$ 
denote by $D_{|z| / 2}$ the disc of radius  $|z|/2$ centered at $0.$
\begin{lem1}\label{LemmaB.2}
Let $\sqrt[\circ]{\cdot }$ be an arbitrary branch of the square root, defined on $z+ D_{|z| / 2},$ $z\in \mathbb{C}\setminus \{0\}.$ 
Then for any  for any $h\in D_{|z| / 2},$ 
\[ 
(i) \quad |\sqrt[\circ]{z+h}-\sqrt[\circ]{z}|\leq \frac{|h|}{|2z|^{1/2}} \quad \text{and}\quad
(ii)\quad  |\sqrt[\circ]{z+h}-\sqrt[\circ]{z}|\leq\sqrt{|h|}/2\, .
\] 
\end{lem1}
\begin{proof}
Note that $\sqrt[\circ]{z+h}-\sqrt[\circ]{z}=\Big(\int_0^1\frac{1}{2}\frac{1}{\sqrt[\circ]{z+th}}\,dt\Big) h$ implying that
\[
|\sqrt[\circ]{z+h}-\sqrt[\circ]{z}|\leq \frac{1}{2}\max_{0\leq t\leq 1}\frac{1}{|z+th|^{\frac{1}{2}}}|h|\leq 
\frac{1}{2}\frac{1}{(|z|/2)^{\frac{1}{2}}}|h|\leq\sqrt{|h|}/2.
\]
\end{proof}

\section{Appendix C: Rough estimates}
In this Appendix, for the convenience of the reader, we prove standard rough asymptotic estimates for the 
Dirichlet eigenvalues $\mu_n$, the periodic eigenvalues $\lambda_n^\pm$, and  of the zeros $\dot\lambda_n$ of 
$\dot\Delta(\lambda, \varphi) $ as $|n|\to \infty$ for potentials $\varphi\in H^1_c$. 
These estimates are needed as a starting point for the proof of our results.

Let $\hat M(x,\lambda):= M(x,\lambda) - E_\lambda(x).$ By \cite{GK}, Theorem 2.3, 
for $0\leq x \leq 1,$ $|\lambda|\geq 1,$ $\varphi\in H^1_c$
\[
|\hat M(x,\lambda)|\leq e^{|\operatorname{Im}\lambda|x}\frac{3\left(1+\|\varphi\|_{L^2}e^{\|\varphi\|_{L^2}}\right)}{|\lambda|} \|
\varphi\|_{H^1}.
\]
Hence for any $\Lambda_{\operatorname{Im}}>0$ and $B>0$ there exists $\Lambda\geq 1$ so that for any $\lambda\in \mathbb{C}$ with 
$|\lambda|\geq \Lambda,$ $|\operatorname{Im}\lambda| \leq \Lambda_{\operatorname{Im}}$

\begin{align}\label{2.1} 
e^{2\Lambda_{\operatorname{Im}}}|\hat M(x,\lambda)|
\leq \frac{1}{2}\frac{\Lambda}{|\lambda|}\qquad
\forall 0\leq x\leq 1, \quad \forall \|\varphi\|_{H^1}\leq B.
\end{align}
It implies that 
\begin{align}\label{2.2}
| M(x,\lambda)| \leq e^{\Lambda_{\operatorname{Im}}}+\frac{1}{2}. 
\end{align}
Using Neumann series to compute $M(x,\lambda)^{-1}$ it then follows that  
\begin{align}\label{2.3}
|M(x,\lambda)^{-1}-E_\lambda(x)^{-1}|\leq 
2e^{2\Lambda_{\operatorname{Im}}} |\hat M(x,\lambda)|\leq \frac{\Lambda}{|\lambda|}
\end{align}
As $E_\lambda(x)^{-1}= E_\lambda(-x)$ this implies that 
\begin{align}
\label{2.4}
| M(x,\lambda)^{-1}| \leq e^{|\operatorname{Im}\lambda|x} +\frac{\Lambda}{|\lambda|} \leq  e^{\Lambda_{\operatorname{Im}}}+1. 
\end{align}
Estimates \eqref{2.1}--\eqref{2.4} will now be used to provide estimates of 
\[\dot M(x,\lambda):= \partial_\lambda M(x,\lambda)\quad \text{and}\quad \ddot  M(x,\lambda):= \partial_\lambda^2 M(x,\lambda)\]
Let $\dot E_\lambda(x)= \partial_\lambda E_\lambda(x)$ and $\ddot E_\lambda(x)= \partial_\lambda^2 E_\lambda(x).$
\begin{lem1}\label{Lemma2.1}
For $\Lambda_{\operatorname{Im}}>0$ and $B>0$ there exist $\Lambda\geq 1$ and $C>0$ so that for any $|\lambda|\geq \Lambda$ with
$|\operatorname{Im}\lambda|\leq\Lambda_{\operatorname{Im}},$ $0\leq x \leq 1,$ and $\|\varphi\|_{H^1}\leq B$
\begin{itemize}
\item[(i)] $|M(x,\lambda)-E_\lambda(x)|\leq \frac{1}{2}\frac{\Lambda}{|\lambda|},\quad |M(x,\lambda)|\leq e^{\Lambda_{\operatorname{Im}}}+\frac{1}{2}$;
\item[(ii)] $|\dot M(x,\lambda)-\dot E_\lambda(x)|\leq C\frac{1}{|\lambda|},\quad |\dot M(x,\lambda)|\leq e^{\Lambda_{\operatorname{Im}}}+C\frac{1}{\Lambda}$;
\item[(iii)] $|\ddot M(x,\lambda)-\ddot E_\lambda(x)|\leq C\frac{1}{|\lambda|},\quad |\ddot M(x,\lambda)|\leq e^{\Lambda_{\operatorname{Im}}}+C\frac{1}{\Lambda}.$
\end{itemize}
\end{lem1}
\begin{proof}
(i) has already been obtained above: see \eqref{2.1}--\eqref{2.2}.
(ii) By \cite{GK}, Corollary 1.5, $\dot M$ satisfies 
\[\dot M(x,\lambda)=-\int_0^x M(x,\lambda)M(t,\lambda)^{-1} RM(t,\lambda)dt\]
Let $(M(x,\lambda)^{-1})\hat{}  := M(x,\lambda)^{-1}- E_\lambda(-x)$ then 
\[
\dot M(x,\lambda)=-\int_0^x \big(E_\lambda(x)+ \hat M(x,\lambda) \big)\big(E_\lambda(-t)+ 
(M(t,\lambda)^{-1})\hat{ } \big) R\big(E_\lambda(t)+ \hat M(t,\lambda) \big)dt.
\]
Note that
\[
-\int_0^x E_\lambda(x)E_\lambda(-t) RE_\lambda(t)dt= -E_\lambda(x)Rx=\dot E_\lambda(x) \, . 
\]
The estimates \eqref{2.1}--\eqref{2.4} then imply that there exists $C>0$ so that 
\[
|\dot M(x,\lambda)-\dot E_\lambda(x)|\leq C\frac{1}{|\lambda|} \qquad
\forall \lambda, \, x, \, \varphi
\]
as in the statement of the lemma. 
As $|\dot E_\lambda(x)|= |E_\lambda(x)||R||x|\leq e^{\Lambda_{\operatorname{Im}}},$ the claimed bound of $\dot M(x,\lambda)$ follows as well.
\\\noindent 
(iii) 
Note that $\ddot M$ satisfies $L\ddot M= \lambda \ddot M+ 2 \dot M.$ Hence by \cite{GK}, Proposition 1.4, 
\[
\begin{array}{l}
\ddot M(x,\lambda)=-\int_0^x M(x,\lambda)M(t,\lambda)^{-1} R\,2\dot M(t,\lambda)dt \\
= -\int_0^x \big(E_\lambda(x)+ \hat M(x,\lambda) \big)\big(E_\lambda(-t)+ (M(t,\lambda)^{-1})\hat{ }  \big) R\,
2\big(\dot E_\lambda(t)+ \hat{\dot M}(t,\lambda) \big)dt
\end{array}
\]
where $\hat{\dot M}(t,\lambda)=\dot M(t,\lambda) -\dot E_\lambda(t).$ As $\dot E_\lambda(t)= - \dot E_\lambda(t)R t$ one sees that 
\[
-\int_0^x E_\lambda(x)E_\lambda(-t) R2\dot E_\lambda(t)dt= \int_0^x E_\lambda(x)R^22t \, dt.
\]
Using that $\ddot E_\lambda(x)= E_\lambda(x)R^2x^2$ one obtains
\[-\int_0^x E_\lambda(x)E_\lambda(-t) R2\dot E_\lambda(t)dt=\ddot E_\lambda(x). \]
The estimates (i),(ii), and \eqref{2.3}--\eqref{2.4} then imply that by choosing $C$ of (ii) larger if needed one gets $|\ddot M(x,\lambda)- \ddot E_\lambda(x)|\leq C\frac{1}{|\lambda|}$ for $\lambda,x$ and $\varphi$ as in the statement of the lemma.
As $|\ddot E_\lambda(x)|= |E_\lambda(x)||R|^2|x|^2\leq e^{\Lambda_{\operatorname{Im}}},$ the claimed bound for $\ddot M(x,\lambda)$ follows as well.
\end{proof}

\noindent
The rough asymptotic estimates for $\mu_n$ and $\dot \lambda_n$ as $|n|\to \infty$ are as follows. 
\begin{lem1}\label{Lemma2.2}
For any  $B>0$ there exists $n_B\geq 1$ so that for any $\varphi \in H^1_c$ with $\|\varphi\|_{H^1}\leq B$
\[
(i)\quad |\mu_n-n\pi|\leq \frac{1}{|n|}\quad \forall |n|\geq n_B \, ;
\qquad
(ii)\quad  |\dot\lambda_n-n\pi|\leq \frac{1}{|n|}\quad \forall |n|\geq n_B \, .
\]
\end{lem1}
\begin{proof}
(i) According to \cite{GK}, Section 5, for any $\varphi\in L^2_c,$ the Dirichlet eigenvalues of $L(\varphi)$ are the zeros (with multiplicities) of 
the characteristic function $\chi_D(\lambda),$ 
\[
\chi_D(\lambda)= 
\frac{1}{2i}\big(m_4(1,\lambda)+m_3(1,\lambda)-m_2(1,\lambda)-m_1(1,\lambda)\big) \, .
\]
Using that by Lemma \ref{Lemma2.1} (i) with $\Lambda_{\operatorname{Im}}=1,$ the absolute value of each entry of 
$M(1,\lambda)-E_\lambda(1)$ is bounded by $\frac{1}{2}\frac{\Lambda}{|\lambda|}$ for any $|\lambda|\geq \Lambda$ with 
$|\operatorname{Im}\lambda|\leq 1$ and $\|\varphi\|_{H^1}\leq B,$ it then follows that 
\[
|2i(\chi_D(\lambda)-\sin \lambda)| =  |(m_4(1,\lambda)-e^{i\lambda})+m_3(1,\lambda)-m_2(1,\lambda)-(m_1(1,\lambda)-e^{-i\lambda})|
\]
can be estimated as $|2i(\chi_D(\lambda)-\sin \lambda)| \leq  \:4\frac{1}{2}\frac{\Lambda}{|\lambda|}$ leading to
\begin{align}
\label{2.5} 
|\chi_D(\lambda)-\sin \lambda|\leq \, \frac{\Lambda}{|\lambda|}.
\end{align}
Increasing $\Lambda$ if needed and arguing in a similar way it follows from Lemma \ref{Lemma2.1} (ii) that 
\begin{align}
\label{2.6} 
|\dot\chi_D(\lambda)-\cos \lambda|\leq 2C\frac{1}{|\lambda|}.
\end{align}
Using that $H^1_c$ embeds compactly into $L^2_c$ one concludes from \cite{GK}, Lemma 5.2, that there exists $n_0\geq 1$ so that
\[
|\mu_n-n\pi| \leq \frac{1}{4}\qquad \forall \, |n|\geq n_0, \; \forall \, \|\varphi\|_{H^1}\leq B.
\]
In particular, $|\operatorname{Im}\mu_n| \leq \frac{1}{4}$ for $|n|\geq n_0.$
Now choose $n_1\geq n_0$ so that
\[
|n_1\pi|-\frac{1}{4}\geq \Lambda
\]
Then $|\mu_n|\geq \Lambda$ for any $|n|\geq n_1$ and hence by \eqref{2.5}--\eqref{2.6} 
\begin{align*}
|\chi_D(\mu_n)-\sin \mu_n| &\leq \, \frac{\Lambda}{|\mu_n|}\leq \, \frac{\Lambda}{|n\pi|-\frac{1}{4}}\\
|\dot\chi_D(\lambda)-\cos \lambda|\leq & \, 2C \frac{1}{|n\pi|-\frac{1}{4}}
 \qquad\forall \; |\lambda-n\pi|\leq \frac{1}{4}.
\end{align*}
Expanding $\chi_D$ at $n\pi$ one gets
\begin{align}
\label{2.7}
\chi_D(\mu_n)-\chi_D(n\pi)= (\mu_n-n\pi)\int_0^1 \dot \chi_D(n\pi +t(\mu_n-n\pi))dt.
\end{align}
Note that $\chi_D(\mu_n)=0$ and by \eqref{2.5}, for any $|n|\geq n_1$
$|\chi_D(n\pi)|\leq \frac{\Lambda}{|n\pi|}$
so that the left hand side of the identity \eqref{2.7} is $O\big(\frac{1}{|n|}\big).$ On the other hand with $x(t)= t(\mu_n-n\pi)$ one has
\[
\cos (n\pi +x(t))=(-1)^n\cos x(t)= (-1)^n\Big(1+ \sum_{k=1}^\infty
\frac{1}{(2k)!}\big(x(t)^2\big)^k\Big)
\]
implying that 
\[
|\cos (n\pi +x(t))- (-1)^n|\leq |x(t)|^2 e^{|x(t)|^2 }\leq \frac{1}{4}\qquad \forall 0\leq t\leq 1
\]
and hence
\[
\Big|\int_0^1 \cos(n\pi+t(\mu_n-n\pi))dt -(-1)^n\Big|\leq \frac{1}{4}\qquad \forall |n|\geq n_1.
\]
By \eqref{2.6} it then follows that 
\[
\Big|\int_0^1 \dot\chi_D(n\pi+t(\mu_n-n\pi))dt -(-1)^n\Big|\leq \frac{1}{4}+ \frac{2C}{|n\pi|-\frac{1}{4}}\,.
\]
Choosing $n_2\geq n_1$ so that
\[\frac{2C}{|n\pi|-\frac{1}{4}}\leq \frac{1}{4}\quad \forall|n|\geq n_2\]
then implies that for any $|n|\geq n_2$
\[
\Big|\int_0^1 \dot\chi_D(n\pi+t(\mu_n-n\pi))\,dt\Big|\geq \frac{1}{2}.
\]
Hence \eqref{2.7} leads to the estimate 
\[|\mu_n-n\pi|\leq 2\frac{\Lambda}{|n\pi|}\qquad \forall |n|\geq n_2,\; \forall \|\varphi\|_{H^1}\leq B,\]
showing that $n_B\geq n_2$ can be chosen as claimed in item (i). 
\\
\noindent 
Concerning (ii) recall that
$\Delta(\lambda) = m_1(1,\lambda)+ m_4(1,\lambda)$. Lemma \ref{Lemma2.1} (ii), (iii), with
$\Lambda_{\operatorname{Im}}=1$ implies that for any $|\lambda|\geq \Lambda$ with 
$|\operatorname{Im}\lambda|\leq 1$ and $\|\varphi\|_{H^1}\leq B,$
\[
|\dot \Delta(\lambda)+ 2\sin \lambda|\leq\frac{2C}{|\lambda|}, \qquad
|\ddot \Delta(\lambda)+ 2\cos \lambda|\leq\frac{2C}{|\lambda|}
\]
where $\ddot\Delta(\lambda)= \partial_\lambda^2 \Delta(\lambda).$ Now argue as in the proof of item (i) to conclude that 
(after choosing $n_B$ larger if needed) (ii) holds. 
\end{proof}
\noindent
As in the case treated in \cite{GK}, Section 6, where $\varphi\in L^2_c,$ the estimates for the periodic eigenvalues are more involved. 
It turns out that the same method of proof as in \cite{GK} works. 
\begin{lem1}\label{Lemma2.3}
For any $B>0$ there exists $n_B\geq 1$ so that for any $\varphi\in H^1_c$
with $\|\varphi\|_{H^1} \leq B$
\[|\lambda_n^\pm-n\pi|\leq \frac{1}{|n|}\qquad \forall |n|\geq n_B.\]
\end{lem1}
\begin{proof}
According to \cite{GK}, Section 6, the periodic eigenvalues of $L(\varphi)$ are the zeroes (with multiplicities) of $\Delta^2(\lambda)-4.$ 
By \cite{GK}, Lemma  6.4, there exists $n_0\geq 1$ so that $|\lambda_n^\pm-n\pi|\leq \frac{1}{8}\;\forall |n|\geq n_0$ and 
$\|\varphi\|_{H^1} \leq B.$ Furthermore, $\Delta(\lambda_n^\pm)= 2 (-1)^{n}\; \forall |n|\geq n_0.$
The cases where $n$ is even and where it is odd are treated in the same way so we concentrate on the even case only. 
One then has $\Delta(\lambda_n^\pm)-2=0.$ Arguing as in the proof of Lemma \ref{Lemma2.2} (ii) one sees that by 
Lemma \ref{Lemma2.1}, applied with $\Lambda_{\operatorname{Im}}=1,$ one has for any $|\lambda|\geq \Lambda$ with 
$|\operatorname{Im}\lambda|\leq 1$ and $\|\varphi\|_{H^1}\leq B,$
\begin{align}\label{2.9}
|\Delta(\lambda)- 2\cos \lambda|\leq \frac{\Lambda}{|\lambda|} \quad \text{and}\quad
|\dot \Delta(\lambda)+ 2\sin \lambda|,
|\ddot \Delta(\lambda)+ 2\cos \lambda|\leq\frac{2C}{|\lambda|}
\end{align}
where $\ddot\Delta(\lambda)= \partial_\lambda^2 \Delta(\lambda).$
By Lemma \ref{Lemma2.2}, one can choose $n_1\geq \max(n_0,8)$ so that $|n_1\pi|\geq 4 \Lambda$ and 
\begin{align}
\label{2.10}
|\lambda_n^\pm-n\pi|,|\mu_n-n\pi|\leq \frac{1}{|n|}\leq \frac{1}{8}\qquad \forall
|n|\geq n_1.
\end{align}
By \eqref{2.9} it then follows that for $|n|\geq n_1,$
\begin{align}
\label{2.11}
|\ddot \Delta(\lambda)+ 2\cos \lambda|\leq\frac{2C}{|\lambda|}\qquad \forall \lambda
\mbox{ \rm{ with } }|\lambda-n\pi | \leq  \frac{1}{8}. 
\end{align} 
Expanding $\Delta$ at $\dot\lambda_n$ up to order two and then evaluating the expansion at $\lambda_n^\pm$ one gets, using that $\dot \Delta(\dot \lambda_n)=0$ and $\Delta(\lambda_n^\pm)=2$
\begin{align}
\label{2.12} 2- \Delta(\dot \lambda_n)= (\lambda_n^\pm-\dot\lambda_n)^2\int_0^1 \ddot \Delta(\dot \lambda_n+t(\lambda_n^\pm-\dot \lambda_n))\cdot (1-t)dt.
\end{align}
In order to use this identity for estimating $|\lambda_n^\pm-\dot\lambda_n|,$ we need to bound the absolute value of the latter integral away from zero. Let $x(t)=\dot\lambda_n-n\pi+ t (\lambda_n^\pm- \dot\lambda_n).$ As $x(t)= (1-t)(\dot\lambda_n-n\pi)+ t(\lambda_n^\pm-n\pi)$ one has $|x(t)|\leq \frac{1}{4}$ for $|n|\geq n_1$ and $0\leq t\leq 1.$
Together with the estimate
\begin{align}
\label{2.14}
|\cos\big(n\pi +x(t)\big)-\cos n\pi|\leq \sum_{k=1}^\infty \frac{1}{(2k)!}\big( |x(t)|^2\big)^k \leq |x(t)|^2e^{|x(t)|^2}
\end{align}
it then follows that for $|n|\geq n_1$ with $n$ even
\[\Big|
\int_0^1  2 \cos (\dot \lambda_n+t(\lambda_n^\pm-\dot \lambda_n))\cdot (1-t)\,dt- 2\int_0^1 (1-t)\,dt
\Big|\leq \frac{1}{4}.\]
Combined with \eqref{2.11} one then gets 
\begin{align}
\label{2.15}
\Big|
\int_0^1 \ddot \Delta(\dot \lambda_n+t(\lambda_n^\pm-\dot \lambda_n))\cdot (1-t)\,dt-1
\Big|
\leq \frac{1}{4}+\frac{2C}{|n\pi|-\frac{1}{8}}
\end{align}
By choosing $n_2\geq n_1$ so that 
\begin{align}\label{2.16}
\frac{2C}{|n\pi|-\frac{1}{8}}\leq \frac{1}{4}\qquad \forall |n|\geq n_2
\end{align} 
one concludes that for $|n|\geq n_2$ with $n$ even and $\|\varphi\|_{H^1}\leq B,$
\[
\Big|\int_0^1 \ddot \Delta(\dot \lambda_n+t(\lambda_n^\pm-\dot \lambda_n))\cdot (1-t)\,dt\Big|\geq \frac{1}{2}
\]
and in turn, by \eqref{2.12}, 
\begin{align}
\label{2.17}
|\lambda_n^\pm-\dot\lambda_n|^2 \leq 2|\Delta(\dot \lambda_n)-2| \, .
\end{align} 
To get the claimed asymptotics of $\lambda_n^\pm$ for $n$ even we need to show that $\Delta(\dot \lambda_n)-2= O\big(\frac{1}{n^2}\big)$. 
First we prove that $\Delta(\mu_n)-2= O\big(\frac{1}{n^2}\big).$
To this end recall from \cite{GK}, Lemma 6.6, that $\Delta^2(\mu_n)-4= \delta^2(\mu_n).$ Write $\Delta^2(\mu_n)-4$ as a product $(\Delta(\mu_n)-2)(\Delta(\mu_n)+2).$ In order to bound $|\Delta(\mu_n)+2|$ away from $0,$ note that 
\[|\Delta(\mu_n)-2|\leq |2\cos\mu_n-2|+ |\Delta(\mu_n)- 2\cos\mu_n|,\]
\[|2\cos\mu_n-2|\leq |\mu_n-n\pi|^2e^{|\mu_n-n\pi|^2}\leq \frac{1}{8}\quad \forall|n|\geq n_2,\; n \,  \text{even}\]
where we used \eqref{2.14}, and by \eqref{2.10}, 
\[
\frac{\Lambda}{|\mu_n|}\leq \frac{\Lambda}{|n\pi|-\frac{1}{8}}\leq \frac{1}{8}\qquad \forall |n|\geq n_2.
\]
We then conclude from \eqref{2.9} that 
\[|\Delta(\mu_n)-2|\leq \frac{1}{8}+ \frac{\Lambda}{|\mu_n|}\leq \frac{1}{2}\qquad \forall|n|\geq n_2,\; n \,  \text{ even }\]
implying that
$|\Delta(\mu_n)+2|\geq 4- |\Delta(\mu_n)-2| \geq 1$ and in turn 
\[
|\Delta(\mu_n)-2|\leq \frac{|\delta(\mu_n)|^2}{|\Delta(\mu_n)+2|}
\leq |\delta(\mu_n)|^2.
\]
By Lemma \ref{Lemma2.1} (i) and Lemma \ref{Lemma2.2} (i)
\[
|\delta(\mu_n)|\leq |m_2(1,\mu_n)|+ |m_3(1,\mu_n)|\leq \frac{\Lambda}{|\mu_n|}\leq \frac{\Lambda}{|n|}\qquad \forall|n|\geq n_2
\]
yielding
\begin{align}\label{2.18}|\Delta(\mu_n)-2|\leq \frac{\Lambda^2}{n^2}\qquad \forall|n|\geq n_2,\; n \,  \text{ even}.
\end{align}
Finally we estimate $\Delta(\dot \lambda_n)-2$ by evaluating at $\mu_n$ the expansion of $\Delta(\lambda)$ in $\dot\lambda_n$ 
\[
\Delta(\mu_n)- \Delta(\dot \lambda_n)= (\mu_n-\dot\lambda_n)^2\int_0^1 \ddot \Delta(\dot \lambda_n+t(\mu_n-\dot \lambda_n))\cdot (1-t)dt.\]
By \eqref{2.10}, $(\mu_n-\dot\lambda_n)^2\leq \frac{4}{n^2}.$ Arguing as in the proof of \eqref{2.15} one sees that 
\[\Big|
\int_0^1 \ddot \Delta(\dot \lambda_n+t(\mu_n-\dot \lambda_n))\cdot (1-t)dt
\Big|\leq 2.\]
Hence $|\Delta(\mu_n)- \Delta(\dot \lambda_n)|\leq \frac{8}{n^2}$ and when combined with \eqref{2.18}, the inequality
\eqref{2.17} yields $|\lambda_n^\pm-\dot\lambda_n|^2 \leq 2 \big(\frac{\Lambda^2}{n^2}+ \frac{8}{n^2}
\big)$ or $|\lambda_n^\pm-\dot\lambda_n| \leq \frac{2\Lambda+4}{|n|}.
$
Using the inequality $|\dot\lambda_n- n\pi|\leq \frac{1}{|n|}$ of \eqref{2.10} and $|\lambda_n^\pm-n\pi|\leq |\dot\lambda_n- n\pi|+ |\lambda_n^\pm-\dot \lambda_n|$
one gets 
\[|\dot\lambda_n- n\pi|\leq \frac{1}{|n|}\big(1+2\Lambda+4\big)\]
yielding the claimed statement for $n$ even.
\end{proof}


\end{document}